\documentclass[reqno]{amsart}

\usepackage{graphicx}
\usepackage{amsfonts}
\usepackage{amssymb}
\usepackage{amsthm}
\usepackage{amsmath}
\usepackage[all]{xy}
\usepackage{multicol}

\theoremstyle{plain}

\newtheorem{theorem}{Theorem}[section]
\newtheorem{lemma}[theorem]{Lemma}
\newtheorem{proposition}[theorem]{Proposition}
\newtheorem{corollary}[theorem]{Corollary}
\newtheorem*{property-P}{Property P}

\theoremstyle{definition}

\newtheorem{definition}[theorem]{Definition}

\newtheorem{remark}[theorem]{Remark}
\newtheorem{notation}[theorem]{Notation}
\newtheorem{example}[theorem]{Example}

\def\comp{\raisebox{0.2mm}{\ensuremath{\scriptstyle\circ}}}

\renewcommand{\mapsto}{\longmapsto}
\renewcommand{\to}{\longrightarrow}
\newcommand{\To}{\Longrightarrow}
\newcommand{\del}{\partial}
\newcommand{\im}{\ensuremath{\mathsf{Im\,}}}
\renewcommand{\Im}{\ensuremath{\mathrm{Im}}}
\newcommand{\coker}{\ensuremath{\mathsf{Coker\,}}}
\newcommand{\Cok}{\ensuremath{\mathrm{Cok}}}
\renewcommand{\ker}{\ensuremath{\mathsf{Ker\,}}}
\newcommand{\Coeq}{\ensuremath{\mathrm{Coeq}}}
\newcommand{\coeq}{\ensuremath{\mathsf{Coeq\,}}}

\newcommand{\Bc}{\ensuremath{\mathcal{B}}}
\newcommand{\Ac}{\ensuremath{\mathcal{A}}}
\newcommand{\Pc}{\ensuremath{\mathcal{P}}}
\newcommand{\Ab}{\ensuremath{\mathsf{Ab}}}
\newcommand{\ab}{\ensuremath{\mathrm{ab}}}
\newcommand{\nil}{\ensuremath{\mathrm{nil}}}
\newcommand{\Nil}{\ensuremath{\mathsf{Nil}}}
\newcommand{\Sol}{\ensuremath{\mathsf{Sol}}}
\newcommand{\sol}{\ensuremath{\mathrm{sol}}}
\newcommand{\xmod}{\ensuremath{\mathrm{xmod}}}
\newcommand{\PXM}{\ensuremath{\mathsf{PXMod}}}
\newcommand{\XM}{\ensuremath{\mathsf{XMod}}}
\newcommand{\Xc}{\ensuremath{\mathcal{X}}}
\newcommand{\Zc}{\ensuremath{\mathcal{Z}}}

\newcommand{\G}{\ensuremath{\mathbb{G}}}
\newcommand{\Gn}{\ensuremath{\mathbb{G}_{n}}}
\newcommand{\Sc}{\ensuremath{\mathcal{S}}}

\newcommand{\Z}{\ensuremath{\mathbb{Z}}}
\newcommand{\N}{\ensuremath{\mathbb{N}}}

\newcommand{\Cat}{\ensuremath{\mathsf{Cat}}}
\newcommand{\Set}{\ensuremath{\mathsf{Set}}}
\newcommand{\Two}{\ensuremath{\mathsf{2}}}
\newcommand{\Gp}{\ensuremath{\mathsf{Gp}}}
\newcommand{\cod}{\ensuremath{\mathrm{cod\,}}}
\newcommand{\dom}{\ensuremath{\mathrm{dom\,}}}
\newcommand{\ext}{\ensuremath{\mathrm{Ext}}}
\newcommand{\Ext}{\ensuremath{\mathsf{Ext}}}
\newcommand{\Extn}{\ensuremath{\mathsf{Ext}^{n}\!}}
\newcommand{\NExt}{\ensuremath{\mathsf{NExt}}}
\newcommand{\Reg}{\ensuremath{\mathsf{Reg}}}
\newcommand{\Arr}{\ensuremath{\mathsf{Arr}}}
\newcommand{\Arrn}{\ensuremath{\mathsf{Arr}^{n}\!}}
\newcommand{\Regn}{\ensuremath{\mathsf{Reg}^{n}\!}}
\newcommand{\CExt}{\ensuremath{\mathsf{CExt}}}

\newcommand{\Hom}{\ensuremath{\mathsf{Hom}}}
\newcommand{\Ch}{\ensuremath{\mathsf{Ch}}}

\newcommand{\Simpa}{\mathcal{SA}}
\newcommand{\Ec}{\ensuremath{\mathcal{E}}}
\newcommand{\set}[1]{\langle #1 \rangle}
\renewcommand{\hom}{\ensuremath{\mathrm{Hom}}}
\newcommand{\op}{\ensuremath{\mathrm{op}}}
\newcommand{\noproof}{\hfill \qed}

\newbox\skewpullbackbox
\setbox\skewpullbackbox=\hbox{\xy 0;<1mm,0mm>: \POS(4,0)\ar@{-} (-4,0) \ar@{-} (8,4)
\endxy}
\def\skewpullback{\copy\skewpullbackbox}

\newbox\ksewpullbackbox
\setbox\ksewpullbackbox=\hbox{\xy 0;<1mm,0mm>: \POS(0,-8)\ar@{-} (0,-4) \ar@{-} (4,-4)
\endxy}

\newbox\pullbackbox
\setbox\pullbackbox=\hbox{\xy 0;<1mm,0mm>: \POS(4,0)\ar@{-} (0,0) \ar@{-} (4,4)
\endxy}
\def\pullback{\copy\pullbackbox}

\newbox\pushoutbox
\setbox\pushoutbox=\hbox{\xy 0;<1mm,0mm>: \POS(0,4)\ar@{-} (0,0) \ar@{-} (4,4)
\endxy}

\makeindex
\def\theindex{\section*{Index of notations}
 \def\item{\par\hangindent 2em}
 \begin{multicols}{3}}
\def\endtheindex{\end{multicols}}

\begin{document}

\newdir{>>}{{}*!/3.5pt/:(1,-.2)@^{>}*!/3.5pt/:(1,+.2)@_{>}*!/7pt/:(1,-.2)@^{>}*!/7pt/:(1,+.2)@_{>}}
\newdir{ >>}{{}*!/8pt/@{|}*!/3.5pt/:(1,-.2)@^{>}*!/3.5pt/:(1,+.2)@_{>}}
\newdir{ |>}{{}*!/-3.5pt/@{|}*!/-8pt/:(1,-.2)@^{>}*!/-8pt/:(1,+.2)@_{>}}
\newdir{ >}{{}*!/-8pt/@{>}}
\newdir{>}{{}*:(1,-.2)@^{>}*:(1,+.2)@_{>}}
\newdir{<}{{}*:(1,+.2)@^{<}*:(1,-.2)@_{<}}

\title{Higher Hopf formulae for homology\\
via Galois Theory}

\author{Tomas Everaert}
\address{Vakgroep Wiskunde\\
Vrije Universiteit Brussel\\
Pleinlaan~2\\
1050~Brussel\\
Belgium}
\email{teveraer@vub.ac.be}

\author{Marino Gran}
\address{Laboratoire de Math\'ematiques Pures et Appliqu\'ees\\
Universit\'e du Littoral C\^ote d'Opale\\
50,~Rue F.~Buisson\\
62228~Calais\\
France}
\email{gran@lmpa.univ-littoral.fr}

\author{Tim Van der Linden}
\address{Vakgroep Wiskunde\\
Vrije Universiteit Brussel\\
Pleinlaan~2\\
1050~Brussel\\
Belgium}
\email{tvdlinde@vub.ac.be}

\thanks{The first author's research is financed by a Ph.D.\ grant of the Institute of Promotion of Innovation through Science and Technology in Flanders (IWT-Vlaanderen). The third author's research is partly financed by the Eduard \v Cech Center for Algebra and Geometry, Brno, Czech Republic (project number LC505). Both would like to thank the LMPA for its kind hospitality during their stay in Calais.}


\begin{abstract}
\noindent We use Janelidze's Categorical Galois Theory to extend Brown and Ellis's higher Hopf formulae for homology of groups to arbitrary semi-abelian monadic categories. Given such a category $\Ac$ and a chosen Birkhoff subcategory $\Bc$ of $\Ac$, thus we describe the Barr-Beck derived functors of the reflector of $\Ac$ onto $\Bc$ in terms of centralization of higher extensions. In case $\Ac$ is the category $\Gp$ of all groups and $\Bc$ is the category $\Ab$ of all abelian groups, this yields a new proof for Brown and Ellis's formulae. We also give explicit formulae in the cases of groups vs.\ $k$-nilpotent groups, groups vs.\ $k$-solvable groups and precrossed modules vs.\ crossed modules. 
\end{abstract}

\maketitle

\setcounter{tocdepth}{1}
\tableofcontents
\pagestyle{myheadings}{\markboth{TOMAS EVERAERT, MARINO GRAN AND TIM VAN DER LINDEN },\markright{HIGHER HOPF FORMULAE FOR HOMOLOGY VIA GALOIS THEORY}}
\section*{Introduction}\label{Section-Introduction}

Generalizing Hopf's formula~\cite{Hopf} for the second integral homology group to higher dimensions is a well-studied problem, that still deserves to be better understood from a categorical perspective. Partial results were originally obtained by Conrad~\cite{Conrad}, Rodicio~\cite{Rodicio} and St\"ohr~\cite{Stohr}, and the first complete solution---a formula describing $H_{n}$ for all $n$---is due to Brown and Ellis~\cite{Brown-Ellis}. Their work was recently extended by Donadze, Inassaridze and Porter in the paper~\cite{Donadze-Inassaridze-Porter}. Whereas Brown and Ellis use topological methods, the latter proof is entirely algebraic, and also considers the case of groups vs.\ $k$-nilpotent groups instead of just groups vs.\ abelian groups.

The aim of our present paper is two-fold: giving a conceptual and elementary proof of the higher Hopf formulae, while at the same time placing them in a very general framework. In our opinion, the simplest approach to a formula for $H_{n}$ is a proof by induction on $n$. Now even for groups, such an approach naturally leads to the use of categorical methods: the familiar category of groups must be left for more general ones. On the other hand, consequent reasoning along such lines gives a lot of added generality for free. 

How general can we go? Just to give an idea: Brown and Ellis's formulae describe $H_{n}$ for groups vs.\ abelian groups; Donadze, Inassaridze and Porter add groups vs.\ $k$-nilpotent groups for arbitrary $k$; we add groups vs.\ $k$-solvable groups, Lie algebras vs.\ abelian Lie algebras, rings vs.\ zero rings, precrossed modules vs.\ crossed modules, etc. Furthermore, even for groups, this yields a proof that's less complicated than the existing ones, and essentially amounts to an application of the Hopf formula for the second homology object and some standard diagram chasing arguments. While based on the same basic idea of using higher dimensional extensions and, in particular, higher presentations of an object, the main difference between our method and previous ones is that we can use an inductive argument, because our formula for $H_{2}$ also holds in categories of higher extensions. 

That such an approach is possible is due to the existence of the appropriate categorical framework, Janelidze, M\'arki and Tholen's \emph{semi-abelian categories}~\cite{Janelidze-Marki-Tholen} and Borceux and Bourn's \emph{homological categories}~\cite{Borceux-Bourn}. These were introduced to capture the fundamental homological properties of the categories of groups, rings, Lie algebras, crossed modules etc.\ much in the same way as abelian categories do for modules over a ring or sheaves of abelian groups. Our work confirms how well the notion of semi-abelian (or homological) category fulfils this promise. 

An important ingredient to understanding the higher Hopf formulae is Janelidze's insight that \emph{centralization of higher extensions} yields the objects that occur in these formulae~\cite{Janelidze:Double, Janelidze:Hopf-talk}. For instance, consider a group $A$ presented by a double extension $f \colon q \to p$
\[
 \vcenter{\xymatrix{ F \ar@{-{ >>}}[d]_-{q} \ar@{-{ >>}}[r] & F/K_2 \ar@{-{ >>}}[d]^-{p}\\
F/K_1 \ar@{-{ >>}}[r] & A}}
\]
where $K_1$ and $K_{2}$ are normal subgroups of $F$ satisfying $A \cong F/K_1 \cdot K_2$ and such that $F$, $F/K_1$ and $F/K_2$ are free groups. Then the third homology group with coefficients in the group of integers $\Z$ is given by
\[
H_{3} (A,\Z) \cong \frac{[F,F]\cap K_1\cap K_2}{L_{2}[f]},
\]
where the object $L_{2}[f]=[K_1\cap K_2, F] \cdot [K_1,K_2]$, a kind of higher commutator, is defined by the \emph{centralization} 
\[
\xymatrix{F/L_{2}[f] \ar@{.{ >>}}[d] \ar@{.{ >>}}[r] & {F/K_2} \ar@{-{ >>}}[d]\\
{F/K_1} \ar@{-{ >>}}[r] & A}
\]
of $f$---the double extension $f$ universally turned into a central one. Of course, such a \emph{double presentation} $f$ always exists: consider the underlying set/free group comonad $\G = (G,\delta ,\epsilon)$ and take the diagram
\[
\vcenter{\xymatrix{ G^{2}A  \ar@{-{ >>}}[d]_-{G\epsilon_{A}} \ar@{-{ >>}}[r]^-{\epsilon_{GA}} & GA \ar@{-{ >>}}[d]^-{\epsilon_{A}}\\
GA \ar@{-{ >>}}[r]_-{\epsilon_{A}} & A}}
\]
for $f$. In the case of groups, Janelidze realized that also the higher Hopf formulae may be interpreted in these terms~\cite{Janelidze:Hopf-talk}. As far as we know, our paper is the first attempt to use this idea for proving the formulae. Next to the concept of semi-abelian categories, part of the fundamental theory is provided by the paper~\cite{EverVdL2}, where a proof along the same lines is given of the Hopf formula for the second homology object, and some of the needed homological tools are developed.

A categorical theory of central extensions was developed by Janelidze and Kelly in~\cite{Janelidze-Kelly} as an application of Janelidze's Categorical Galois Theory~\cite{Janelidze:Pure}. This theory is modelled on the situation where $\Ac$ is a variety of universal algebras and $\Bc$ a given subvariety of $\Ac$, and allows one to classify the extensions in $\Ac$ that are \emph{central} with respect to this chosen subvariety $\Bc$. The idea of relative centrality---which goes back to the work of the Fr\"ohlich school, see e.g.,~\cite{Froehlich, Lue, Furtado-Coelho}---has for leading example the variety of groups with its subvariety of abelian groups. Janelidze and Kelly's theory is general enough to include the case where $\Ac$ is any semi-abelian category and $\Bc$ is any given Birkhoff subcategory of $\Ac$. Modelling the notion of centralization of higher extensions in a categorical way, we are forced to extend this theory, but still within the framework of Categorical Galois Theory. The resulting process of centralization of higher extensions then provides us with the higher commutators that occur in the higher Hopf formulae. 

In the last section we explain how these commutators may be calculated explicitly, and we do so in some specific cases: groups vs.\ abelian groups, groups vs.\ $k$-nilpotent groups, groups vs.\ $k$-solvable groups, precrossed modules vs.\ crossed modules. In the important example of precrossed modules the homology objects are described by the same formula (see Theorem~\ref{Thm-pxmod}): indeed, the so-called Peiffer commutator plays the same role, in the category of precrossed modules, as the usual commutator of normal subgroups does in the category of groups.

\section{Semi-abelian and homological categories}\label{Semi-abelian categories}

As pointed out above, semi-abelian and homological categories were introduced to capture the homological properties of those categories ``sufficiently close'' to the category $\Gp$ of all groups. In this section we briefly recall their definition and basic properties.

It is important to note that in general, the difference between a semi-abelian category and an abelian one is quite vast: in an abelian category, every morphism may be factored as a cokernel followed by a kernel; any hom-set $\hom (B,A)$ carries an abelian group structure; binary products and binary coproducts coincide. None of these properties holds true for the category $\Gp$ of groups: the first one because not every subgroup is a normal subgroup, the second one essentially because the pointwise product of two group homomorphisms need not be a homomorphism, and the third one because a group $A$ with $A\times A\cong A+A$ is always trivial. In view of these differences, it is easily understood why the definition of semi-abelian category might sound unfamiliar at first. Nevertheless, the link with the notion of semi-abelian category is simple and precise: a category $\Ac$ is abelian if and only if both $\Ac$ and its dual category ${\Ac}^{\op}$ are semi-abelian~\cite{Janelidze-Marki-Tholen}.

\begin{definition}\label{Definition-Semi-Abelian}
A category $\Ac$ is \emph{semi-abelian} when it is pointed, Barr exact and Bourn protomodular and has binary coproducts~\cite{Janelidze-Marki-Tholen}. $\Ac$ is \emph{homological} ~\cite{Borceux-Bourn} when it is pointed, regular and Bourn protomodular.
\end{definition}

Of course some explanation is needed. First of all, in a semi-abelian category, all finite limits and colimits exist. In particular, there is a terminal object $1$ and an initial object $0$, and it is possible to construct finite products and coproducts, equalizers and coequalizers, pullbacks and pushouts. $\Ac$ being \emph{pointed} means that $0\cong 1$, i.e., there is a \emph{zero} object: an object that is both initial and terminal. A map is called \emph{zero} when it factors over $0$; given any two objects $A$ and $B$, there is a unique zero map from $B$ to $A$. This makes it possible to consider kernels and cokernels: given a morphism $f\colon {B\to A}$ in $\Ac$, a \emph{kernel} $\ker f\colon {K[f]\to B}$ of $f$ is a pullback of $0\to A$ along $f$ and dually, a \emph{cokernel} $\coker f\colon {A\to \Cok [f]}$ is a pushout of $B\to 0$ along $f$. 

$\Ac$ being \emph{Barr exact} means that it is regular, and such that every internal equivalence relation in $\Ac$ is a kernel pair~\cite{Barr}. We start by commenting on the regularity, and later come back to the other condition. Recall that a morphism is called a \emph{regular epimorphism} when it is a coequalizer of some pair of arrows. Having finite limits and coequalizers of kernel pairs, $\Ac$ is \emph{regular} when moreover the regular epimorphisms of $\Ac$ are pullback-stable. In a regular category, image factorizations exist: any morphism $f\colon {B\to A}$ can be factored as a regular epimorphism ${B\to \Im[f]}$ followed by a monomorphism $\im f\colon {\Im [f]\to A}$ called the \emph{image} of~$f$; this factorization is unique up to isomorphism. A morphism $f$ such that $\im f$ is a kernel is called \emph{proper}.

In this pointed and regular context, $\Ac$ is \emph{Bourn protomodular}~\cite{Bourn1991} if and only if the (regular) \emph{Short Five Lemma} holds: this means that for every commutative diagram
\[
\xymatrix{K[f_{0}] \ar@{{ >}->}[r]^-{\ker f_{0}} \ar[d]_-k & B_{0} \ar@{->>}[r]^-{f_{0}} \ar[d]^-b & A_{0} \ar[d]^-a \\ 
K[f] \ar@{{ >}->}[r]_-{\ker f} & B \ar@{->>}[r]_-{f} & A}
\]
such that $f$ and $f_{0}$ are regular epimorphisms, $k$ and $a$ being isomorphisms entails that $b$ is an isomorphism. This implies that every regular epimorphism is in fact a cokernel (of its kernel). Accordingly, we can define exact sequences as follows. A~sequence of morphisms $(f_{i})_{i\in I}$
\[
\xymatrix{\dots \ar[r] & A_{i+1} \ar[r]^-{f_{i+1}} & A_{i} \ar[r]^-{f_{i}} & A_{i-1}\ar[r] & \dots}
\]
is called \emph{exact at $A_{i}$} if $\im f_{i+1}=\ker f_{i}$. (In particular, then $f_{i+1}$ is proper.) It is called \emph{exact} when it is exact at $A_{i}$, for all $i\in I$. Sequence~\ref{Short-Exact-Sequence} below 
\begin{equation}\label{Short-Exact-Sequence}
\xymatrix{
0 \ar[r] & K \ar@{{ >}->}[r]^-{k} & B  \ar@{-{>>}}[r]^-{f} & A \ar[r] & 0}
\end{equation}
is exact if and only if it represents $(k,f)$ as a \emph{short exact sequence}: $k=\ker f$ and $f=\coker k$.

Finally, under the above assumptions, $\Ac$ will be Barr exact if and only if the direct image of a kernel along a regular epimorphism is a kernel: given any kernel $k$ and any regular epimorphism $f$ in $\Ac$, their composition $f\comp k$ is proper---its image $\im f\comp k$ is a kernel.

Examples of semi-abelian categories include all abelian categories; any variety of $\Omega$-groups: amongst the operations defining it, there is a group operation and a unique constant (the unit of the group operation), in particular, the categories of groups, (non-unital) rings, (pre)crossed modules, Lie algebras over a field, commutative algebras; compact Hausdorff groups; $C^{*}$-algebras; the dual of the category of pointed sets.

For some applications, having a homological category instead of a semi-abelian one will be sufficient. For instance, the category of topological groups and of torsion-free abelian groups are homological but not semi-abelian, and we shall meet another example in Section~\ref{Section-Higher-Extensions} below. Yet, already in this context the basic homological diagram lemma's hold: the \emph{Snake Lemma}, the \emph{$3\times 3$-Lemma}, etc. \cite{Bourn2001}

\section{Galois structures}\label{Section-Galois-Structures}
In this section we recall the basic definition of \emph{categorical Galois structure}, which is crucial for the study of higher central extensions. In particular we give a useful sufficient condition for a Galois structure to be \emph{admissible} in the sense of this theory. We refer the reader to  the monograph~\cite{Borceux-Janelidze} by Borceux and Janelidze, and in particular to its introduction, for the historical background that led to the development of the theory, as well as for the details of several interesting examples of admissible Galois structures.

\begin{definition}~\cite{Janelidze:Precategories}
A \emph{Galois structure} \index{Gamma@$\Gamma$}$\Gamma=(\Ac,\Bc,\Ec,\Zc,I,H)$ consists of two categories, $\Ac$ and $\Bc$, an adjunction
\[
\xymatrix{\Ac \ar@<1 ex>[r]^-{I} \ar@{}[r]|-{} & \Bc, \ar@<1 ex>[l]^-H \ar@{}[l]|{\perp}}
\] 
and classes $\Ec$ and $\Zc$ of morphisms of $\Ac$ and $\Bc$ respectively, such that:
\begin{enumerate}
\item $\Ac$ has pullbacks along arrows in $\Ec$;
\item 
$\Ec$ and $\Zc$ contain all isomorphisms, are closed under composition and are pullback-stable;
\item
$I(\Ec) \subset \Zc$;
\item
$H(\Zc) \subset \Ec$;
\item
the counit $\epsilon$ is an isomorphism;
\item
each $A$-component $\eta_{A}$ of the unit $\eta$ belongs to $\Ec$.
\end{enumerate}
An element of $\Ec$ is called an \emph{extension}.
\end{definition}

\begin{example}\label{ExampleGroups}
Let $\Gp$ be the category of groups and $\Ab$ its full reflective subcategory of abelian groups:
\[
\xymatrix{\Gp \ar@<1 ex>[r]^-{\ab } \ar@{}[r]|-{} & {\Ab.} \ar@<1 ex>[l]^-{H} \ar@{}[l]|{\perp}}
\] 
This adjunction determines a Galois structure $\Gamma=(\Gp,\Ab,\Ec,\Zc,\ab,H)$ where $H$ is the inclusion functor, $\ab$ is the abelianization functor, and $\Ec$ and $\Zc$ are the classes of surjective homomorphisms in $\Gp$ and in $\Ab$, respectively. For any group $A$, the $A$-component of the unit of the adjunction is given by the canonical quotient $\eta_A \colon A \to \ab (A)={A}/{[A,A]}$. In particular every $\eta_A$ belongs to $\Ec$.
\end{example}

For an object $A$ of $\Ac$, let us denote $\Ec(A)$ the full subcategory of the slice category $\Ac/A$ determined by the arrows $B\to A$ in $\Ec$.  If $a \colon A_{0}\to A$ is an arrow, we will write $a^* \colon \Ec (A)\to \Ec (A_{0})$ for the functor that sends an extension $f\colon B\to A$ to its pullback $a^{*}f$ along $a$. Since maps in $\Ec$ are pullback-stable, $a^*$ is well-defined. 

For every object $A$ of $\Ac$, the adjunction $(I,H)$ gives rise to an adjunction
\[
\xymatrix{\Ec(A) \ar@<1 ex>[r]^-{I^{A}} \ar@{}[r]|-{} & \Zc(IA), \ar@<1 ex>[l]^-{H^{A}} \ar@{}[l]|-{\perp}}
\] 
with $I^A$ defined by $I^{A}f=If$ on arrows $f$ of $\Ec(A)$, and $H^{A}$ defined by $H^{A}x=(\eta_{A})^{*}Hx$ on arrows $x$ of $\Zc(IA)$. The unit and counit of this adjunction will be denoted by $\eta^A$ and $\epsilon^A$, respectively.  

\begin{definition}
A Galois structure $\Gamma$ is called \emph{admissible} if, for every object $A$ of~$\Ac$, the counit $\epsilon^{A}\colon  I^{A} H^{A} \To 1_{\Zc(IA)}$ is an isomorphism.
\begin{equation}\label{Admissiblediagram}
\vcenter{\xymatrix@1{
A \times_{HIA} HX \ar`u[r]`[rrr]^{\pi_2}[rrr]  \ar[rr]^-{\eta_{A \times_{HIA} HX}} \ar[d]_{\pi_1=H^{A}x=\eta_{A}^{*}Hx} && HI(A \times_{HIA} HX) \ar@{.>}[r]^-{\epsilon_{x}^{A}} \ar[d]_{HI\pi_1} & HX \ar[d]^{Hx} \\
A \ar[rr]_{\eta_A} && HIA \ar@{=}[r] & HIA}}
\end{equation}
\end{definition}

\begin{example}
The Galois structure $\Gamma=(\Gp,\Ab,\Ec,\Zc,\ab,H)$ of Example~\ref{ExampleGroups} is admissible~\cite{Janelidze:Pure}. More generally, any subvariety of a congruence modular variety determines an admissible Galois structure where $\Ec$ and $\Zc$ are the classes of surjective homomorphisms~\cite{Janelidze-Kelly}.
\end{example}

\begin{definition}\label{Definition-Admissible}
Let $\Ec$ be a class of pullback-stable regular epimorphisms in a category $\Ac$. A full replete $\Ec$-reflective subcategory $\Bc$ of $\Ac$
\[
\xymatrix{\Ac \ar@<1 ex>[r]^-{I} \ar@{}[r]|-{} & \Bc \ar@<1 ex>[l]^-H \ar@{}[l]|-{\perp}}
\] 
is a \emph{strongly $\Ec$-Birkhoff} subcategory of $\Ac$ if, for every $f \colon B \to A$ in $\Ec$, all arrows in the diagram 
\[
\vcenter{\xymatrix{B \ar@/^/[drr]^-{\eta_B} \ar@/_/[drd]_-{f} \ar@{.>}[rd]|{(f,\eta_B)} \\
& A \times_{HIA} HIB \ar@{}[rd]|<<{\pullback} \ar[r] \ar[d] & HIB \ar[d]^-{HIf} \\
& A \ar[r]_-{\eta_A} & HIA}}
\]
are in $\Ec$. In particular, the outer diagram is a \emph{regular pushout} in the sense of Bourn (see~\cite{Bourn2003, Carboni-Kelly-Pedicchio} and Diagram~\ref{Diagram-Regular-Pushout} below).
\end{definition}

When $\Ac$ is a regular category and $\Ec$ is the class of regular epimorphisms in~$\Ac$, the notion of strongly $\Ec$-Birkhoff subcategory of $\Ac$ is stronger than the classical Birkhoff property. This property requires $\Bc$ to be closed in $\Ac$ under subobjects and regular quotients or, equivalently, the outer diagram above to be a pushout of regular epimorphisms. A Birkhoff subcategory (in the classical sense) of a variety of universal algebras is the same thing as a subvariety. The notions of Birkhoff subcategory and of strongly $\Ec$-Birkhoff subcategory coincide as soon as $\Ac$ is exact and satisfies the Mal'tsev property~\cite{Carboni-Kelly-Pedicchio}: every internal reflexive relation in $\Ac$ is an equivalence relation. For instance, if $\Ac$ is a Mal'tsev variety, then the subvarieties of $\Ac$ are strongly (regular epi)-Birkhoff subcategories. Recall that every semi-abelian category is exact Mal'tsev~\cite{Bourn1996, Borceux-Bourn}.

In general, a Galois structure with the property that the adjunction satisfies the strongly $\Ec$-Birkhoff property is always admissible:
\begin{proposition}
Let $\Gamma=(\Ac,\Bc,\Ec,\Zc,I,H)$ be a Galois structure, where $\Bc$ is a strongly $\Ec$-Birkhoff subcategory of $\Ac$, with $\Ec$ a given class of pullback-stable regular epimorphisms. Then $\Gamma$ is an admissible Galois structure.
\end{proposition}
\begin{proof}
This is a consequence of the following two facts. On the one hand, we have that $\eta^A_f$ is an epimorphism for all $f\colon B\to A$ in $\Ec(A)$, by the strongly $\Ec$-Birkhoff property of $\Bc$. On the other hand, $H^A$ reflects isomorphisms since, again by the strongly $\Ec$-Birkhoff property, $\eta_A$ is a pullback-stable regular epimorphism: see, for instance, Proposition~1.6 in~\cite{Janelidze-Sobral-Tholen}. To see that these facts indeed imply that $\epsilon^{A}_{x}$ is an isomorphism, for any $x\colon X\to IA$ in $\Zc(IA)$, consider the triangular identity $(H^A\epsilon^A_x) \comp \eta^A_{H^A x}=1_{H^A x}$. Now, $\eta^A_{H^A x}$ is both a split monomorphism and an epimorphism, hence an iso. This implies that $H^A\epsilon^A_x$ is an iso and, finally, $\epsilon^A_x$ is an iso, since $H^A$ reflects isomorphisms. 
\end{proof}

\section{Higher extensions}\label{Section-Higher-Extensions}

We now restrict our attention to the situation where $\Ac$ is a homological category. We shall be considering higher-dimensional arrows in $\Ac$; they are the objects of the following categories.

\begin{definition}
Let $\Two$ be the category generated by a single map $\varnothing\to \{1 \}$. For any $n$, write $\Two^{n}$ for the $n$-fold product $\Two \times \cdots \times \Two$, and denote the functor category $\Hom(\Two^n,\Ac)$ as $\Arrn\Ac$\index{ArrnA@$\Arrn\Ac$}.
\end{definition}

Recall that a \emph{regular epimorphism} is a coequalizer of some pair of arrows. As a first approach to higher-dimensional extensions, we could use higher-dimensional regular epimorphisms:

\begin{definition} 
Let $\Reg \Ac$ be the full subcategory of the category of arrows in $\Ac$ whose objects are the regular epimorphisms. Denote $\Reg^{0}\!\Ac =\Ac$ and $\Regn\Ac=\Reg(\Reg^{n-1}\!\Ac)$ for $n\geq 1$. By an \emph{$n$-fold regular epimorphism} we mean an object of $\Regn\Ac$\index{RegnA@$\Regn\Ac$}, i.e., a regular epimorphism between $(n-1)$-fold regular epimorphisms. 
\end{definition}
For brevity, we shall say \emph{$n$-regular epimorphism} instead of $n$-fold regular epimorphism.

Note that a double (=~$2$-) regular epimorphism of $\Ac$ may be considered as a commutative square in $\Ac$ and, in general, an $n$-regular epi as a particular kind of commutative $n$-dimensional diagram. This is a formal consequence of the fact that the functor $(\cdot)\times \Two\colon \Cat\to \Cat$ is left adjoint to the functor $\Hom (\Two,\cdot)\colon \Cat \to\Cat$: thus an object of $\Reg^{2}\!\Ac$, being a functor ${\Two \to \Hom (\Two ,\Ac)}$, corresponds to a functor ${\Two\times \Two\to \Ac}$, i.e., a commutative square in $\Ac$; by induction, an object of $\Regn\Ac$ can be seen as a functor $\Two^{n}\to \Ac$. 

It is not difficult to show that an object of $\Reg^{2}\!\Ac$ is a pushout square in $\Ac$, and in general, an object of $\Regn\Ac$ is an $n$-dimensional cube in $\Ac$ of which all (two-dimensional) faces are pushouts.

When $\Ac$ is a semi-abelian category, the category $\Arrn\Ac$ is of course semi-abelian, as is any category of $\Ac$-valued presheaves. On the other hand, while $\Regn\Ac$ is still homological (see~\cite{Tomasthesis}), it is no longer semi-abelian. Indeed, it is well known that for a (non-trivial) abelian category $\Ac$ the category $\Reg\Ac$ is \emph{not} abelian; hence it cannot be exact, since it is obviously additive. To see this, consider an object $A$ of $\Ac$ such that $A\neq 0$ and write $\tau_A$ for the unique arrow $A\to 0$. Then the diagram 
\[
\xymatrix{
A \ar@{=}[r]   \ar@{=}[d] & A  \ar[d]\\
A \ar[r] & 0 }
\]    
represents an epimorphism $1_A\to \tau_A$ in $\Reg\Ac$ which is not a normal epimorphism (a cokernel). It follows that $\Reg\Ac$ is not exact whenever $\Ac$ is a non-trivial abelian category. It is then easily seen  that also for $n\geq 2$, $\Regn\Ac$ need not be semi-abelian.

It is well known that when, in a regular category, a commutative square of regular epimorphisms
\begin{equation}\label{Diagram-Double-RegEpi}
\vcenter{\xymatrix{B_{0} \ar@{->>}[r]^-{f_{0}} \ar@{->>}[d]_-{b} & A_{0} \ar@{->>}[d]^-a \\
B \ar@{->>}[r]_-f & A}}
\end{equation}
is a pullback, it is a pushout. In a regular category, a commutative square of regular epimorphisms is called a \emph{regular pushout} when the comparison map $r\colon B_{0}\to P$ to a pullback
\begin{equation}\label{Diagram-Regular-Pushout}
\vcenter{\xymatrix{B_{0} \ar@/^/@{->>}[drr]^-{f_{0}} \ar@/_/@{->>}[drd]_-{b} \ar@{.>}[rd]|r \\
& P \ar@{}[rd]|<{\pullback} \ar@{.>>}[r] \ar@{.>>}[d] & A_{0} \ar@{->>}[d]^-a \\
& B \ar@{->>}[r]_-f & A}}
\end{equation}
of $f$ along $a$ is a regular epimorphism. An important aspect of these regular pushouts is made clear by the next result.

\begin{proposition}\cite{Bourn2003}\label{Proposition-Rotlemma-Kernels}
Consider, in a homological category, a commutative diagram of exact sequences, such that $f$, $f_{0}$, $b$ and $a$ are regular epimorphisms:
\[
\xymatrix{0 \ar[r] & K[f_{0}] \ar[d]_-{k} \ar@{{ >}->}[r] & B_{0} \ar@{-{>>}}[d]_-{b} \ar@{-{>>}}[r]^-{f_{0}} & A_{0} \ar@{-{>>}}[d]^-{a} \ar[r] & 0\\
0 \ar[r] & K[f] \ar@{{ >}->}[r] & B \ar@{-{>>}}[r]_-{f} & A \ar[r] & 0.}
\]
The right hand square is a regular pushout if and only if $k$ is regular epi.\noproof
\end{proposition}

For a regular epimorphism of regular epimorphisms that is a regular pushout in $\Ac$, this means that its kernel is computed degreewise or, equivalently, that degreewise taking kernels in the diagram above induces a morphism $k$ in $\Ac$ that is an object of $\Reg\Ac$. 
In~\cite{Carboni-Kelly-Pedicchio} Carboni, Kelly and Pedicchio show that a regular category $\Ac$ is Barr exact and Mal'tsev if and only if in $\Ac$, every pushout of regular epimorphisms is a regular pushout. In particular, a semi-abelian category has this property: a pushout of two regular epimorphisms always exists, and it is a regular pushout. On the other hand, the failure of $\Regn\Ac$ to be exact implies that an $(n+2)$-regular epimorphism is \emph{not} the same as a regular pushout in $\Regn\Ac$, for $n\ge 1$. This difference gives rise to the notion of $n$-extension.

\begin{definition}\label{Definition-extension}
Let $\Ac$ be a homological category. A \emph{$0$-fold extension in~$\Ac$} is an object of $\Ac$ and a \emph{$1$-fold extension} is a regular epimorphism in $\Ac$. An \emph{$n$-fold extension} is an object $(f_0,f)$ of $\Arrn\Ac$ such that all arrows in the induced diagram~\ref{Diagram-Regular-Pushout} are $(n-1)$-fold extensions. The $n$-fold extensions determine a full subcategory $\Extn\Ac$\index{ExtnA@$\Extn\Ac$} of $\Arrn\Ac$.
\end{definition}
Again, we shall say \emph{$n$-extension} instead of $n$-fold extension.

As an important technical result, we shall need Proposition~\ref{Proposition-Higher-Rotlemma}, which generalizes Proposition~\ref{Proposition-Rotlemma-Kernels} to higher extensions. But first we show that some constructions in $\Extn \Ac$ may be performed in $\Arrn \Ac$, which will prove very useful later on.

\begin{proposition}\label{Proposition-Higher-Ext-Pullback-Stable}
Let $\Ec^n$\index{En@$\Ec^{n}$} be the class of $(n+1)$-extensions in $\Ac$. Then:
\begin{enumerate}
\item pullbacks in $\Extn \Ac$ along maps in $\Ec^n$ exist, are pullbacks in $\Arrn\Ac$ and, in particular if $n\geq 1$, are computed degreewise in $\Ext^{n-1}\! \Ac$;
\item maps in $\Ec^n$ are stable under pulling back in $\Extn \Ac$;
\item $\Ec^n$ is closed under composition.
\end{enumerate}
\end{proposition}
\begin{proof}
If $n=0$ the statements are true because $\Ac$ is regular. Now consider $n\geq 1$ and suppose that the statements hold for $k<n$. 

To prove the first two statements, consider the following pullback in $\Arrn\Ac$ of an $n$-extension $(\alpha_{0} ,\alpha)$ along an arbitrary morphism $(f_{0},f)$. 
\[
\xymatrix@=.5cm{& B'_{0} \ar@{}[rd]|<{\skewpullback} \ar[dd]|<<<{\hole}|-{\hole}_(.75){b'} \ar@{->}[rr]^-{f'_{0}} \ar@{-{ >>}}[ld]_-{\beta_{0}} && A'_{0} \ar@{-{ >>}}[dd]^-{a'} \ar@{-{ >>}}[ld]_-{\alpha_{0}}\\
B_{0} \ar@{-{ >>}}[dd]_-{b} \ar@{->}[rr]^(.75){f_{0}} && A_{0} \ar@{-{ >>}}[dd]^(.25){a}\\
& B' \ar@{->}[rr]^(.25){f'}|(.5){\hole} \ar@{}[rd]|<{\skewpullback} \ar@{-{ >>}}[ld]^-{\beta} && {A'} \ar@{-{ >>}}[ld]^-{\alpha}\\
B \ar@{->}[rr]_-{f} && A}
\]
Using that $(n-1)$-extensions are closed under pulling back, we see that $\beta$ and $\beta_{0}$ are $(n-1)$-extensions. We are to show that $(\beta_{0} ,\beta)$ is an $n$-extension; then in particular, it is a pullback in $\Extn \Ac$ of $(\alpha_{0} ,\alpha)$ along $(f_{0},f)$.

Pulling back $\alpha$ along $a$, $\beta$ along $b$ yields the commutative diagram
\[
\xymatrix{ B'_{0} \ar@{->}[r]^-{f'_{0}} \ar@{->}[d]_-{(\beta_{0},b')} & A'_{0} \ar@{-{ >>}}[d]^-{(\alpha_{0},a')}\\
{B_{0}\times_{B}B'} \ar@{->}[r]_-{f_{0}\times_{f}f'} & {A_{0}\times_{A}A'.}}
\]
By assumption, $(\alpha_{0},a')$ is an $(n-1)$-extension; moreover, this diagram is a pullback, hence, by the induction hypothesis, $(\beta_{0},b')$ is an $(n-1)$-extension. Now $b'$, as the composite of $(\beta_{0},b')$ with the second projection of the pullback of $b$ along $\beta$, is a composite of $(n-1)$-extensions, hence, by the induction hypothesis, is an $(n-1)$-extension itself.

To prove the third statement, considering two composable $n$-extensions 
\[
\xymatrix{C_{0} \ar@{-{ >>}}[r]^-{g_{0}} \ar@{-{ >>}}[d]_-{c} & B_{0} \ar@{-{ >>}}[d]_-b \ar@{-{ >>}}[r]^-{f_{0}} & A_{0} \ar@{-{ >>}}[d]^-a \\
C \ar@{-{ >>}}[r]_-g & B \ar@{-{ >>}}[r]_-f & A}
\]
we immediately see that also $f_{0}\comp g_{0}$ and $f\comp g$ are $(n-1)$-extensions. Since $(n-1)$-extensions are stable under pulling back, we moreover get that every arrow in the diagram 
\[
\xymatrix{C_{0} \ar@{-{ >>}}[r]^-{g_{0}} \ar@{-{ >>}}[d]_-{r_{g}} & B_{0} \ar@{=}[d] \\
P_{g} \ar@{}[rd]|<{\pullback} \ar@{-{ >>}}[r]^-{\overline{g}} \ar@{-{ >>}}[d]_-{\overline{r_f}} & B_{0} \ar@{-{ >>}}[r]^-{f_{0}} \ar@{-{ >>}}[d]_-{r_{f}} & A_{0} \ar@{=}[d] \\
P \ar@{}[rd]|<<{\pullback} \ar@{-{ >>}}[d]_-{\overline{\overline{a}}} \ar@{-{ >>}}[r]^-{} & P_{f} \ar@{}[rd]|<{\pullback} \ar@{-{ >>}}[d]_-{\overline{a}} \ar@{-{ >>}}[r]^-{\overline{f}} & A_{0} \ar@{-{ >>}}[d]^-a \\
C \ar@{-{ >>}}[r]_-g & B \ar@{-{ >>}}[r]_-f & A}
\]
is an $(n-1)$-extension, which proves our claim.
\end{proof}

In particular, the kernel $K[f]$ in $\Arrn \Ac$ of an $(n+1)$-extension $f\colon {B\to A}$ is always an $n$-extension. We are going to show in Proposition \ref{Proposition-Higher-Rotlemma} that the converse is also true. Suppose a map $f$ in $\Extn \Ac$ has a kernel (in $\Arrn \Ac$) that is an $n$-extension; suppose moreover that this map $f$ is regular epi in $\Arrn \Ac$: then this map $f$ is in $\Ec^{n}$. 

\begin{lemma}\label{Lemma-g}
Suppose that $n\geq 0$. Consider, in $\Arrn\Ac$, a commutative diagram of exact sequences
\[
\xymatrix{0 \ar[r] & K[f_{0}] \ar[d]_-{k} \ar@{{ >}->}[r] & B_{0} \ar[d]_-{b} \ar@{-{>>}}[r]^-{f_{0}} & A \ar@{=}[d] \ar[r] & 0\\
0 \ar[r] & K[f] \ar@{{ >}->}[r] & B  \ar@{-{>>}}[r]_-{f} & A \ar[r] & 0.}
\]
If $f$ and $f_{0}$ are $(n+1)$-extensions, then $k\in \Ec^{n}$ if and only if $b\in\Ec^{n}$.
\end{lemma}
\begin{proof}
Whenever, in the above diagram, $b$, $f$ and $f_{0}$ are $(n+1)$-extensions, the right hand side square is an $(n+2)$-extension, hence $k$ is an $(n+1)$-extension by Proposition~\ref{Proposition-Higher-Ext-Pullback-Stable}. On the other hand, by Proposition~8 in~\cite{Bourn2001}, if $k$ is a regular epimorphism then so is $b$. By induction, it is then easy to verify the following: if $B$ is an $n$-extension, and both $f_0$ and $k$ are $(n+1)$-extensions, then also $b$ is in~$\Ec^{n}$.
\end{proof}

\begin{lemma}\label{Lemma-Higher-Rotlemma}
Suppose that $n\geq 0$. Consider, in $\Arrn\Ac$, a commutative diagram of exact sequences, such that $f$, $f_{0}$, $b$ and $a$ are $(n+1)$-extensions:
\[
\xymatrix{0 \ar[r] & K[f_{0}] \ar[d]_-{k} \ar@{{ |>}->}[r] & B_{0} \ar@{-{ >>}}[d]_-{b} \ar@{-{ >>}}[r]^-{f_{0}} & A_{0} \ar@{-{ >>}}[d]^-{a} \ar[r] & 0\\
0 \ar[r] & K[f] \ar@{{ |>}->}[r] & B  \ar@{-{ >>}}[r]_-{f} & A \ar[r] & 0.  }
\]
The right hand square is an $(n+2)$-extension if and only if $k$ is an $(n+1)$-extension. In particular, if $f$ is an $(n+2)$-extension, then its kernel $K[f]$ in $\Ext^{n+1}\!\Ac$ exists and is computed degreewise.
\end{lemma}
\begin{proof}
It follows from the previous Lemma.
\end{proof}

\begin{lemma}\label{Lemma-Higher-Ext-Compo}
Consider a pair of composable arrows $f\colon B\to A$ and $g\colon C\to B$ in $\Arrn\Ac$. If $f \comp g$ is an $(n+1)$-extension and $B$ an $n$-extension, then $f$ is an $(n+1)$-extension.\noproof
\end{lemma}

\begin{proposition}\label{Proposition-Higher-Rotlemma}
Suppose that $n\geq 0$. Consider, in $\Arrn\Ac$, a short exact sequence 
\begin{equation}
\xymatrix{
0 \ar[r] & K \ar@{{ >}->}[r]^-{k} & B  \ar@{-{>>}}[r]^-{f} & A \ar[r] & 0}
\end{equation}
such that $B$ is an $n$-extension. Then $f$ is an $(n+1)$-extension if and only if $K$ is an $n$-extension.\noproof 
\end{proposition}

This result shows that for $n\geq 1$, an $n$-extension is the same thing as a sequence~\ref{Short-Exact-Sequence} in $\Ext^{n-1}\!\Ac$ exact in $\Arr^{n-1}\!\Ac$. A kernel in $\Arr^{n-1}\!\Ac$ between $(n-1)$-extensions always has a cokernel, and this cokernel is an $n$-extension. However, in general, an exact sequence in $\Ext^{n-1}\!\Ac$ does not determine an extension! Thus it may also be shown that an $n$-dimensional arrow in $\Ac$ is an $n$-extension exactly when it is an \emph{exact $n$-presentation} in the sense of~\cite{Donadze-Inassaridze-Porter}. Moreover, this allows one to describe extensions in terms of kernels alone. 

\begin{corollary}\label{Corollary-Higher-Denormalized-Rotlemma}
Suppose that $n\geq 0$. Consider, in $\Arrn\Ac$, an exact fork
\[
\xymatrix{R[f] \ar@<.5ex>[r]^-{\pi_{1}} \ar@<-.5ex>[r]_-{\pi_{2}} & B \ar@{-{>>}}[r]_-{f} & A}
\]
such that $B$ is an $n$-extension. Then $f$ is an $(n+1)$-extension if and only if $R[f]$ is an $n$-extension. 
\end{corollary}
\begin{proof}
This is a consequence of Proposition~\ref{Proposition-Higher-Rotlemma} and~\ref{Proposition-Higher-Ext-Pullback-Stable}: it suffices to note that $\ker f =\pi_{2}\comp \ker \pi_{1}$.
\end{proof}

\begin{example}\label{Example-Split-Epi-extension}
A split epimorphism in $\Extn\Ac$ is an $(n+1)$-extension. 
\end{example}

\begin{notation}\label{Notation-Zeroes}
The forgetful functor $\Upsilon \colon{\Ext\Ac \to \Ac}$ that maps an extension $f\colon {B\to A}$ to the object $B$ has an obvious right adjoint, namely the functor $\Psi\colon {\Ac \to \Ext \Ac}$ that sends an object $B$ of $\Ac$ to the extension $B\to 0$. Composition yields adjunctions\index{Psin@$\Psi^{n}$}
\[
\xymatrix{\Extn\Ac \ar@<1ex>[r]^-{\Upsilon^{n} } \ar@{}[r]|-{} & \Ac. \ar@<1ex>[l]^-{ \Psi^{n}} \ar@{}[l]|-{\perp}}
\]
\end{notation}

\section{Higher central extensions}

We shall now establish a sequence of Galois structures $\Gamma_{n}$, such that each determines the next one in the following way: $\Gamma_{0}$ is induced by a Birkhoff subcategory $\Bc$ of a semi-abelian category $\Ac$, and $\Gamma_{n}$ is a structure on $\Extn \Ac$ with class of extensions $\Ec^{n}$ and in which the adjunction is given by centralization of $n$-extensions with respect to $\Gamma_{n-1}$. 

\subsection{Trivial, central and normal extensions}\label{Subsection-Trivial-Central-Normal}
Let $\Gamma=(\Ac,\Bc,\Ec,\Zc,I,H)$ be a Galois structure such that $\Bc$ is an $\Ec$-Birkhoff subcategory of $\Ac$. From now on, we will omit the inclusion $H$ from our notations (and write $\subseteq$ for the functor $H$). We shall adopt the terminology of~\cite{Janelidze-Kelly}, calling \emph{trivial}, \emph{central} and \emph{normal} the following types of extensions:
\begin{definition}\label{Definition-Central-Extension}
Let $f\colon {B\to A}$ be an extension. One says that $f$ is
\begin{enumerate}
\item \emph{a trivial extension} (with respect to $\Gamma$), when the next square is a pullback;
\[
\xymatrix{
B \ar[r]^-{\eta_{B}}  \ar[d]_-{f} & IB \ar[d]^-{If} \\
A \ar[r]_-{\eta_{A}} & IA}
\]
\item a \emph{central extension}, when there exists an $a\colon{A_{0}\to A}$ in $\Ec(A)$ such that $a^{*}f$ is trivial;
\item a \emph{normal extension}, when the first projection $\pi_1 \colon R[f] \to B$ (or, equivalently, the second projection $\pi_2$) is a trivial extension.
\end{enumerate}
\end{definition}

It is clear that every normal extension is central. $\Gamma$ being admissible implies that, moreover, every trivial extension is normal. To see this, note that $f$ being trivial implies that $R[f]\cong R[If]\times_{IB} B$; but $R[If]$ is in $\Bc$, because $\Bc$ is a full replete reflective subcategory of $\Ac$; hence the admissibility of $\Gamma$ entails $R[If]\cong IR[f]$. 

\begin{example}\label{centralextensionsforgroupsxample}
Let us again consider Example~\ref{ExampleGroups} and its Galois structure $\Gamma=(\Gp,\Ab,\Ec,\Zc,\ab,\subseteq )$. One easily sees that the trivial extensions are exactly the surjective homomorphisms of groups $f \colon B \to A$ with the property that the restriction $\hat{f} \colon [B,B] \to [A,A]$ of $f$ to the derived subgroups is an isomorphism. It was shown in~\cite{Janelidze:Pure} that an extension $f \colon B \to A$ that is central with respect to~$\Gamma$ is the same thing as a central extension in the classical sense: its kernel $K[f]$ is contained in the centre $Z(B)$ of $B$.  
\end{example}

\subsection{The functors $I_{n}$ and $J_{n}$}\label{Subsection-In-Jn}\index{In@$I_{n}$}\index{Jn@$J_{n}$}
Let $\Ac$ be a semi-abelian category, $\Bc$ a Birkhoff subcategory of $\Ac$ and $I$ the reflector $\Ac\to \Bc$, with unit $\eta$. Of course, we may also consider $I$ as a functor ${\Ac\to\Ac}$. Putting $J(A)=K[\eta_A]$ for each object $A$ of~$\Ac$ yields another functor $J\colon{\Ac\to\Ac}$ as well as a short exact sequence of functors
\[
\xymatrix{
0 \ar[r] & J \ar@{{ >}->}[r]^-{\mu} & 1_{\Ac} \ar@{-{>>}}[r]^-{\eta} & I \ar[r] & 0.}
\]
 From this sequence, we construct a short exact sequence of functors ${\Arr\Ac\to\Arr\Ac}$ as follows: let $f$ be an extension $B\to A$ and $(\pi_1,\pi_2)$ its kernel pair. Put $J_1[f]= K[J\pi_1]$ and $J_1f= \Psi J_1[f]$.
\[
\xymatrix{
J_1[f]=K[J\pi_1] \ar@{}[dr]|{\texttt{(i)}} \ar@{{ >}->}[r]^-{\ker J\pi_{1}} \ar[d] & JR[f] \ar[d]^{\mu_{R[f]}} \ar@<0.5 ex>[r]^-{J\pi_1} \ar@<-0.5 ex>[r]_-{J\pi_2} & JB \ar[d]_{\mu_B}\\
K[\pi_1] \ar@{{ >}->}[r]_-{\ker \pi_{1}} & R[f] \ar@<0.5 ex>[r]^-{\pi_1} \ar@<-0.5 ex>[r]_-{\pi_2} & B}
\]
Clearly, this defines a functor $J_{1}\colon {\Arr\Ac\to\Arr\Ac}$. Furthermore, let us define $\mu^1_f=(\mu_B\comp J\pi_2\comp\ker J\pi_1, \alpha_A)$, where $\alpha_A$ is the unique arrow from $0$ to $A$, so that $\eta^1_f= (\rho^1_f,1_A)=\coker \mu^1_f$ yields a short exact sequence 
\[
\xymatrix{
0 \ar[r] & J_1 \ar@{{ >}->}[r]^-{\mu^1} & 1_{\Arr\Ac} \ar@{-{>>}}[r]^-{\eta^1} & I_1 \ar[r] & 0}
\]
of functors ${\Arr\Ac\to\Arr\Ac}$. Indeed, $\mu^1_f$ is a monomorphism because both $\mu_B$ and $J\pi_2\comp\ker J\pi_1$ are monomorphisms: $\mu_B$ by assumption, and $J\pi_2\comp\ker J\pi_1$ because it is the normalization of the reflexive, hence effective equivalence relation $(JR[f], J\pi_1,J\pi_2)$. Since $\mu_B$ is a monomorphism, the square \texttt{(i)} is a pullback, hence $\mu_{R[f]}\comp \ker J\pi_1$ is a normal monomorphism, as an intersection of normal monomorphisms. It follows that $\pi_2\comp\mu_{R[f]}\comp \ker J\pi_1$, the regular image of $\mu_{R[f]}\comp \ker J\pi_1$ along $\pi_2$, is normal in $B$; hence so is $\mu^1_f$.  

Since $\Arr\Ac$ is semi-abelian as soon as $\Ac$ is, we may repeat this process inductively in order to obtain, for each $n\ge 0$, a short exact sequence\index{etan@$\eta^{n}$}\index{mun@$\mu^{n}$} 
\[
\xymatrix{
0 \ar[r] & J_n \ar@{{ >}->}[r]^-{\mu^n} & 1_{\Arrn\Ac} \ar@{-{>>}}[r]^-{\eta^n} & I_n \ar[r] & 0}
\]
of functors $\Arrn\Ac\to\Arrn\Ac$. Here we put $J_0=J$ and $I_0=I$. As in the case $n=1$, we write $J_n[f]$ for the domain of $J_nf$ and $I_n[f]=B/J_n[f]$ for the domain of $I_nf$, for any $n$-extension $f\colon {B\to A}$. Also, we define $\rho^n_f\colon {B\to I_n[f]}$ via $\eta_f^n=(\rho^n_f, 1_A)$.\index{rhon@$\rho^{n}$}\index{Inn@$I_{n}[\cdot]$}\index{Jnn@$J_{n}[\cdot]$} Hence, for each $n$-extension $f\colon {B\to A}$, we have, in $\Arr^{n-1}\!\Ac$, a short exact sequence
\[
\xymatrix{
0 \ar[r] & J_n[f] \ar@{{ >}->}[r]^-{\ker \rho_f^n} & B  \ar@{-{>>}}[r]^-{\rho_f^n} & I_n[f] \ar[r] & 0.}
\]

\subsection{The Galois structures $\Gamma_{n}$}\label{Subsection-Structures-Gamma_{n}}\index{Gamman@$\Gamma_{n}$}
Given a Birkhoff subcategory $\Bc$ of a semi-abelian category $\Ac$, we denote $\Gamma_0=(\Ac,\Bc,\Ec,\Zc,I,\subseteq )$ the associated Galois structure: $\Ec$ and $\Zc$ are the classes of regular epimorphisms in $\Ac$ and in $\Bc$, respectively. It is well known that for this structure, the central and normal extensions coincide~\cite{Janelidze-Kelly}. 

Let $\CExt_{\Bc}\Ac=\CExt_{\Bc}^1\Ac$\index{CExtBA@$\CExt_{\Bc}\Ac$} denote the full subcategory of $\Ext^1\!\Ac$ whose objects are the $\Gamma_0$-central extensions. Write  $\Ec^{1} = \Ext^2\!\Ac$. We are going to show that $I_1$ (co)restricts to a reflector $\Ext\Ac\to\CExt_{\Bc}\Ac$ and that $\CExt_{\Bc}\Ac$ is a strongly $\Ec^1$-Birkhoff subcategory of $\Ext\Ac$. In particular, this gives rise to a Galois structure 
\[
\Gamma_1=(\Ext\Ac, \CExt_{\Bc}\Ac, \Ec^1,\Zc^1, I_1, \subseteq _1),
\]
where $\Zc^1$ consists of all $2$-extensions in $\CExt_{\Bc}\Ac$ and $\subseteq _1$ denotes the inclusion of $\CExt_{\Bc}\Ac$ into $\Ext \Ac$. Inductively, this may be extended to higher extensions: if 
\[
\Gamma_{n-1}=(\Ext^{n-1}\!\Ac,\CExt^{n-1}_{\Bc}\!\Ac,\Ec^{n-1},\Zc^{n-1},I_{n-1},\subseteq _{n-1})
\]
is the $(n-1)$-th Galois structure in the sequence, write 
\[
\CExt_{\Bc}^n\Ac = \NExt_{\CExt_{\Bc}^{n-1}\!\Ac}(\Ext^{n-1}\!\Ac),
\]
where $\NExt_{\CExt_{\Bc}^{n-1}\!\Ac}$\index{NExtnBA@$\NExt^{n}_{\Bc}\Ac$} is the full subcategory of $\Ext^{n-1}\!\Ac$ of normal extensions with respect to $\Gamma_{n-1}$. Let $\Gamma_{n}$\index{Gamman@$\Gamma_{n}$} be the structure 
\[
(\Extn\Ac,\CExt^{n}_{\Bc}\Ac,\Ec^{n},\Zc^{n},I_{n},\subseteq _{n}),
\]
where $\Ec^{n}$\index{En@$\Ec^{n}$}\index{Zn@$\Zc^{n}$} and $\Zc^{n}$ are all $(n+1)$-extensions in $\Extn \Ac$ and $\CExt_{\Bc}^{n}\Ac$,\index{CExtnBA@$\CExt^{n}_{\Bc}\Ac$} respectively, and\index{subn@$\subseteq_{n}$} $\subseteq_{n}$ denotes the inclusion $\CExt_{\Bc}^{n}\Ac\to \Extn \Ac$. The notation $\CExt_{\CExt_{\Bc}^{n-1}\!\Ac}(\Ext^{n-1} \Ac)$ is explained by Proposition~\ref{Central-is-normal}: for all $n\ge 2$, the normal and central extensions with respect to $\Gamma_{n-1}$ coincide. 

\begin{lemma}
Let $\Bc$ be a Birkhoff subcategory of a semi-abelian category $\Ac$. For each $n\ge 1$, $\CExt_{\Bc}^n\Ac$ is a strongly $\Ec^n$-Birkhoff subcategory of $\Extn\Ac$. The reflector is given by the restriction of $I_n$ to a functor $\Extn\Ac\to\CExt_{\Bc}^n\Ac$.
\end{lemma}
\begin{proof}
If we define $\CExt_{\Bc}^0\Ac=\Bc$, then the case $n=0$ is true by assumption. Let us then suppose that the Lemma holds for all $0\le k\le n-1$. We prove it for $k=n$. 

First of all, $I_n\colon {\Arrn\Ac\to\Arrn\Ac}$ must restrict to a functor ${\Extn\Ac\to \CExt_{\Bc}^n\Ac}$. Note that  $\CExt_{\Bc}^n\Ac$ is well defined by the induction hypothesis. Let $f\colon {B\to A}$ be an $n$-extension. By Lemma~\ref{Lemma-Higher-Ext-Compo}, $I_nf$ is an $n$-extension as well. We must show that it is normal with respect to $\Gamma_{n-1}$. 

Let $(\pi_1,\pi_2)$ denote the kernel pair of $f$ and $(\pi_1',\pi_2')$ the kernel pair of $I_nf$, computed in $\Ext^{n-1}\Ac$, or, equivalently, in $\Arr^{n-1}\Ac$. Consider the following diagram.
\[
\xymatrix{
0 \ar[r] & J_{n-1}R[I_nf] \ar@{{ >}->}[r] \ar[d]_{J_{n-1}\pi_1'}  & R[I_nf] \ar@{}[rd]|{\texttt{(ii)}} \ar@{-{>>}}[r]^-{\eta^{n-1}_{R[I_nf]}}   \ar[d]_{\pi_1'} & I_{n-1}R[I_nf] \ar[r]  \ar[d]^{I_{n-1}\pi_1'}   & 0\\
0 \ar[r] & J_{n-1}I_n[f] \ar@{{ >}->}[r] & I_n[f] \ar@{-{>>}}[r]_-{\eta^{n-1}_{I_n[f]}} & I_{n-1}I_n[f] \ar[r] & 0}
\]
We must prove that the square \texttt{(ii)} is a pullback in $\Ext^{n-1}\Ac$. By Proposition~\ref{Proposition-Higher-Ext-Pullback-Stable}, this is equivalent to proving that it is a pullback in $\Arr^{n-1}\Ac$, taking into account that $\CExt^{n-1}\Ac$ is strongly $\Ec^{n-1}$-Birkhoff in $\Ext^{n-1}\Ac$. We are going to show that $J_{n-1}\pi_1'=J_{n-1}\pi_2'$. Since $J_{n-1}\pi_1'$ and $J_{n-1}\pi_2'$ are jointly monic, this implies that $J_{n-1}\pi_1'$ is a monomorphism hence an iso, so that by Theorem~2.3 in~\cite{Bourn-Janelidze:Semidirect}, \texttt{(ii)} is a pullback. 

Since $\CExt_{\Bc}^{n-1}\Ac$ is a strongly $\Ec^{n-1}$-Birkhoff subcategory of $\Ext^{n-1}\Ac$, $J_{n-1}$ preserves $n$-extensions, hence the left hand downward pointing arrow in the diagram   
\[
\xymatrix{
J_{n-1}R[f] \ar[d] \ar@<0.5 ex>[r]^{J_{n-1}\pi_1} \ar@<-0.5 ex>[r]_{J_{n-1}\pi_2} & J_{n-1}B \ar[d] \ar[r]^-{\mu^{n-1}_B} & B \ar[d]^{\rho^n_f} \\
J_{n-1}R[I_n f]  \ar@<0.5 ex>[r]^{J_{n-1}\pi'_1} \ar@<-0.5 ex>[r]_{J_{n-1}\pi'_2} & J_{n-1}I_n[f] \ar[r]_-{\mu^{n-1}_{I_n[f]}} & I_n[f].}
\]
is an $n$-extension; in particular, it is an epimorphism. Furthermore, $\mu^{n-1}_{I_n[f]}$ is a monomorphism. Hence to show that $J_{n-1}\pi'_1=J_{n-1}\pi'_2$, it suffices to prove 
\[
\rho^n_f\comp \mu^{n-1}_B\comp J_{n-1}\pi_1=\rho^n_f\comp \mu^{n-1}_B\comp J_{n-1}\pi_2. 
\]
But this follows from the fact that $\rho^n_f\comp \mu^{n-1}_B$ factors over $\coeq (J_{n-1}\pi_1,J_{n-1}\pi_2)$: 
\[
\xymatrix{
0 \ar[r] &J_n[f] \ar@{{ >}->}[r] \ar@{=}[d] & J_{n-1}B \ar[d]_{\mu^{n-1}_B} \ar@{-{>>}}[r] & \frac{J_{n-1}B}{J_n[f]} \ar@{.>}[d] \ar[r] & 0 \\
0 \ar[r] &J_n[f] \ar@{{ >}->}[r] & B \ar@{-{>>}}[r]_-{\rho^n_f} & {I_n[f]= \frac{B}{J_n[f]}} \ar[r] & 0.}
\]
Let us then show that $I_n\colon {\Extn\Ac\to \CExt_{\Bc}^n\Ac}$ is a left adjoint. Let $h\colon {f\to g}$ be an arrow in $\Extn\Ac$ with $g$ in $\CExt_{\Bc}^n\Ac$. Let $(\pi_1,\pi_2)$ be the kernel pair of $g$. The normality of $g$ implies that $J_{n-1} \pi_1$ is an isomorphism. It follows that $J_n[g]=K[J_{n-1}\pi_1]=0$, hence $I_ng=g$. Consequently, $I_nh$ gives a factorization $I_nf\to I_ng=g$:
\[
\xymatrix{
f\ar[r]^-{\eta^n_f} \ar[dr]_h & I_nf \ar@{.>}[d]^{I_nh}\\
& g,}
\]
which is unique, because $\eta^n_f$ is an epimorphism.

To see that $\CExt_{\Bc}^n\Ac$ is a strongly $\Ec^n$-Birkhoff subcategory of $\Extn\Ac$, it remains to be shown that $J_{n}$ preserves $(n+1)$-extensions: then the strongly $\Ec^n$-Birkhoff property follows from Proposition~\ref{Proposition-Higher-Rotlemma}. Let $(f_0,f)\colon {b\to a}$ be an $(n+1)$-extension. Denote the kernel pair of $f_0\colon {B_0\to A_0}$ by $(\pi_1,\pi_2)$ and the kernel pair of  $f\colon {B\to A}$ by $(\pi'_1,\pi'_2)$. Since $\CExt_{\Bc}^{n-1}\!\Ac$ is a strongly $\Ec^{n-1}$-Birkhoff subcategory of $\Ext^{n-1}\!\Ac$, $J_{n-1}$ preserves $n$-extensions. Therefore, both $J_{n-1}R[(f_0,f)]$ and $J_{n-1}b$ are $n$-extensions. It follows that the right hand square in the diagram
\[
\xymatrix{
J_n[f_0] \ar[r] \ar[d] & J_{n-1}R[f_0] \ar@{-{ >>}}[d]_-{J_{n-1}R[(f_0,f)]} \ar@<.5ex>[r]^-{J_{n-1}\pi_1} & J_{n-1}B_0 \ar@{-{ >>}}[d]^{J_{n-1}b} \ar@<.5ex>[l] \\
J_n[f] \ar[r] & J_{n-1}R[f]  \ar@<.5ex>[r]^-{J_{n-1}\pi'_1} & J_{n-1}B \ar@<.5ex>[l]}
\]
is a split epi between $n$-extensions, therefore it is an $(n+1)$-extension. We can conclude that the left hand downward pointing arrow is an $n$-extension, hence $J_n(f_0,f)$ is an $(n+1)$-extension.
\end{proof}

We are now going to prove that the objects of the category $\CExt_{\Bc}^n\Ac$ are indeed the $n$-fold central extensions with respect to the Galois structure $\Gamma_{n-1}$. To prove this, we need the next lemma.

\begin{lemma}\label{split-pullback}
For any $n\geq 0$, the reflector $I_n \colon \Extn\Ac \to \CExt_{\Bc}^n\Ac$ preserves pullbacks of split epimorphisms along morphisms in ${\Ec}^{n}$.
\end{lemma}
\begin{proof}
By the Short Five Lemma, an $(n+2)$-extension $(f_{0},f)$, considered as a square~\ref{Diagram-Double-RegEpi} in $\Extn \Ac$, is a pullback if and only if its kernel $K[(f_{0},f)]$ is an isomorphism in $\Arrn\Ac$. Now let $f$ be a split epimorphism in $\Extn \Ac$, and $a$ an $(n+1)$-extension, in $\Extn \Ac$ and let $(f_0,f) \colon b \to a$ denote the extension obtained by pulling back $f$ along $a$. 

Consider the diagram with exact rows
\[
\xymatrix{0 \ar[r] & K[(f_{0},f)] \ar@{{ |>}->}[r] \ar@{-{ >>}}[d] & b \ar@{}[rd]|{\texttt{(iii)}} \ar@{-{ >>}}[r]^-{(f_{0},f)} \ar@{-{ >>}}[d] & a \ar@{-{ >>}}[d] \ar[r] & 0\\
0 \ar[r] & K[I_{n}(f_{0},f)] \ar@{{ |>}->}[r] & I_{n}b \ar@{-{ >>}}[r]_-{I_{n}(f_{0},f)} & I_{n}a \ar[r] & 0}
\]
in $\Arr^{n+1}\!\Ac$. By the strongly $\Ec^{n}$-Birkhoff property of $I_{n}$, the arrows ${b\to I_{n}b}$ and ${a\to I_{n}a}$ are $(n+2)$-extensions. Hence the square \texttt{(iii)} is an $(n+3)$-extension, since $(f_{0},f)$ is a split epimorphism. Since an $\Ec^{n+1}$-quotient of an isomorphism in $\Extn \Ac$ is an isomorphism---isomorphisms are stable under pushing out---we get that $K[I_{n}(f_{0},f)]$ is iso, and hence the $(n+2)$-extension $I_{n}(f_{0},f)= (I_{n}f_{0},I_{n}f)$, considered as a square in $\Arrn\Ac$, is a pullback. 
\end{proof}

\begin{proposition}\label{Central-is-normal}
For any $n\geq 1$, let $\Gamma_n$ be the Galois structure constructed above.
\begin{enumerate}
\item The central extensions with respect to $\Gamma_{n}$ are the normal extensions. Hence the objects in $\CExt_{\Bc}^{n+1}\!\Ac$ are indeed exactly all extensions that are central with respect to $\Gamma_n$.
\item Every central extension that is a split epimorphism is trivial.
\end{enumerate}
\end{proposition}
\begin{proof}
Let $f \colon B \to A$ be a central extension. Since $f^{*}f$ is a split epimorphism, (1) will follow from (2) if we show that central extensions are pullback-stable in $\Extn \Ac$. Now this latter property is an easy consequence of the fact that $\CExt^{n+1}_{\Bc}\!\Ac$ is a strongly $\Ec^{n+1}$-Birkhoff subcategory of $\Ext^{n+1}\!\Ac$, which implies that trivial (and hence also central) extensions are pullback-stable.

Suppose that $f \colon B \to A$ is a split epic central extension. Its centrality means that there exists an arrow $a \colon A_{0} \to A$ in $\Ec^{n+1}$ for which $a^*f$ is trivial. Using this fact and the previous lemma, one concludes that the exterior rectangle and the left hand square in the diagram
\[
\xymatrix{ A_{0} \times_A B \ar@{}[rd]|<<{\pullback} \ar@{-{ >>}}[r]^-{\pi_2} \ar@{-{ >>}}[d]_-{\pi_1} & B \ar@{-{ >>}}[d]_-f \ar@{-{ >>}}[r]^-{\eta^n_B} & I_n B \ar@{-{ >>}}[d]^-{I_n f} \\
A_{0} \ar@{-{ >>}}[r]_-a & A \ar@{-{ >>}}[r]_-{\eta^n_A} & I_n A}
\]
are pullbacks. Since $a$ is a pullback-stable regular epimorphism, it follows that the right hand square is a pullback, and $f$ is a trivial extension.
\end{proof}

\begin{theorem}\label{Theorem-Galois}
If $\Ac$ is a semi-abelian category and $\Bc$ is a Birkhoff subcategory of $\Ac$, then for every $n\geq 1$, the $n$-extensions and central $n$-extensions give rise to a Galois structure $\Gamma_{n}$. This structure is admissible and strongly $\Ec^{n}$-Birkhoff, and provides the corresponding notion of central $(n+1)$-extensions.\noproof  
\end{theorem}

\begin{remark}\label{Remark-Effective-for-Descent}
Like the categories $\Regn \Ac$, the categories $\Extn \Ac$ need not be exact when $n\geq 1$. However, using Proposition~\ref{Proposition-Higher-Ext-Pullback-Stable} and Corollary~\ref{Corollary-Higher-Denormalized-Rotlemma} one easily deduces from the exactness of $\Arrn \Ac$ that every $(n+1)$-extension is an effective descent morphism in $\Extn \Ac$. (For more details on descent theory, we refer the reader to~\cite{Janelidze-Sobral-Tholen}.) Thus the Fundamental Theorem of Categorical Galois Theory from~\cite{Janelidze:Precategories} yields, for any given $n$-extension $a\colon {A_{0}\to A}$, a classification of the $n$-extensions $f$ with codomain $A$ and $a^{*}f$ trivial. Suppose $\Ac$ has enough (regular epi)-projective objects; then any $\Extn \Ac$ has enough $\Ec^{n}$-projectives: for every object $A$ in $\Extn \Ac$, an $(n+1)$-extension $p\colon {P\to A}$ can be chosen such that $P$ is $\Ec^{n}$-projective. Then, since trivial $(n+1)$-extensions are pullback-stable in $\Extn \Ac$, for all $(n+1)$-extensions $f$ and $a$ with codomain $A$, $p^{*}f$ is trivial as soon as $a^{*}f$ is. In this case, the Fundamental Galois Theorem yields a classification of the central $(n+1)$-extensions with a fixed codomain $A$.
\end{remark}

\begin{example}\label{Centralization}
Let $I_1=\ab_1$ be the reflection $\Ext^1\Gp \to \CExt_{\Ab}^1 \Gp$ sending an extension $f \colon B \to A$ of groups to its centralization $\ab_1 f \colon {B}/{[K[f],B]}$ $ \to A$ (see Example~\ref{centralextensionsforgroupsxample}). This reflection gives rise to the Galois structure 
\[
\Gamma_1= (\Ext^{1}\Gp,\CExt^{1}_{\Ab}\Gp,\Ec^{1},\Zc^{1},\ab_{1},\subseteq _{1}).
\]
It was shown by Janelidze~\cite{Janelidze:Double} that the double extensions that are central with respect to this Galois structure $\Gamma_1$ are precisely those extensions~\ref{Diagram-Double-RegEpi} with the property that $[K[f_0],K[b]]= \{ 1 \}$ and $[K[f_0] \cap K[b], B_0]= \{1 \}$. Recently in~\cite{Gran-Rossi} a similar characterization was obtained, valid in the context of Mal'tsev varieties.
\end{example}

\begin{notation}\label{Notation-W}
Let $f\colon B\to A$ be an $n$-extension. Recall from the definition of the $J_n$ that $J_{n}f=\Psi J_{n}[f]$. It is easily seen that $J_{n}f$ lies, moreover, in the image of $\Psi^n$. We write $L_n[f]$\index{Lnn@$L_{n}[\cdot]$} for the object defined via $\Psi^n L_{n}[f]= J_n f$. It follows that only the ``top map'' of the unit  $\eta^{n}_{f}$ at $f$ is not an iso. We denote it by $\eta_{f}$\index{eta@$\eta$}.
\end{notation}

\section{Simplicial extensions}\label{Section-Simplicial-extensions}

We recall the semi-abelian definition of Barr-Beck cotriple homology, and we describe how a suitable comonad on $\Ac$ produces canonical simplicial resolutions of higher extensions in $\Ac$.

\subsection{Homology in semi-abelian categories}\label{Subsection-Homology-Semiabelian}
In a pointed category, a \emph{chain complex} $C={(d_{n}\colon C_{n}\to C_{n-1})}_{n\in \Z}$ is a sequence of maps satisfying $d_{n}\comp d_{n+1}=0$. As in the abelian case, if it exists, the \emph{$n$-th homology object} of $C$ is\index{Hn@$H_{n}$} 
\[
H_{n}C=\Cok[C_{n+1}\to K[d_n]],
\]
the cokernel of the factorization of $d_{n+1}$ over the kernel of $d_{n}$. Now for this definition to be of any interest, one demands that the complex $C$ and the underlying category satisfy some additional properties. Recall that a morphism in a pointed and regular category is called \emph{proper} when its image is a kernel. As soon as the ambient category is, moreover, protomodular, homology of \textit{proper} chain complexes---those with boundary operators of which the image is a kernel---is well-behaved: it characterizes exactness of complexes, and any short exact sequence of proper chain complexes induces a long exact homology sequence~\cite{EverVdL2}. 

This notion of homology may be extended to simplicial objects: one considers the \emph{normalization functor} $N\colon \Simpa \to \Ch\Ac$ that maps a simplicial object $S$ in a pointed category with pullbacks $\Ac$ to its \textit{Moore complex} $NS$, the chain complex with $N_{n}S=0$ for $n<0$, $N_{0}S=S_{0}$,
\[
N_{n} S=\bigcap_{i=0}^{n-1}K[\del_{i}\colon S_{n}\to S_{n-1}]
\]
and boundary operators $d_{n}=\del_{n}\comp \bigcap_{i}\ker \del_{i}\colon N_{n} S\to N_{n-1} S$, for $n\geq 1$. Then $H_{n}S=H_{n}N S$. As shown in~\cite[Theorem~3.6]{EverVdL2}, it is easily seen that when $\Ac$ is semi-abelian, every simplicial object in $\Ac$ has a proper Moore complex.

Let $\Ac$ be an arbitrary category and\index{G@$\G$}\index{epsilon@$\epsilon$} 
\[
\G =(G\colon \Ac \to \Ac ,\quad \delta\colon G\To G^{2} ,\quad \epsilon\colon G\To 1_{\Ac})
\]
a comonad on $\Ac$. Recall the axioms of comonad: $\epsilon_{GA}\comp \delta_{A}=G\epsilon_{A}\comp \delta_{A}=1_{GA}$ and $\delta_{GA}\comp \delta_{A}=G\delta_{A}\comp \delta_{A}$, for any object $A$ of $\Ac$. Putting
\[
\del_{i}=G^{i}\epsilon_{G^{n-i}A}\colon G^{n+1}A\to G^nA, \qquad \sigma_{i}=G^{i}\delta_{G^{n-i}A}\colon G^{n+1}A\to G^{n+2}A ,
\]
for $0\leq i \leq n$, gives the sequence $(G^{n+1}A)_{n\in \N}$ the structure of a simplicial object $\G A$ of~$\Ac$. It has an augmentation $\epsilon_{A}\colon {GA\to A}$; the augmented simplicial object $\epsilon_{A} \colon {\G A\to A}$ is called the \emph{canonical $\G$-simplicial resolution of $A$}. The following naturally generalizes Barr-Beck cotriple homology~\cite{Barr-Beck} to the semi-abelian context.

\begin{definition}\cite{EverVdL2}\label{Definition-Homology}
Let $\Ac$ be a category equipped with a comonad $\G$ and $\Bc$ a semi-abelian category. Let $I\colon \Ac\to \Bc$ be a functor. For $n\geq 1$, the object
\[
H_{n} (A,I)_{\G}=H_{n-1} N I \G A
\]
is the \emph{$n$-th homology object of $A$ (with coefficients in $I$) relative to the cotriple~$\G$}.\index{HnIG@$H_{n}(\cdot,I)_\G$} This defines a functor $H_{n}(\cdot,I)_\G \colon \Ac \to \Bc $, for every $n\geq 1$. 
\end{definition}

\begin{example}\label{Example-Comonads} (cf. Example~\ref{ExampleGroups})
When $\Ac$ is the category of groups, $\ab\colon {\Gp \to \Ab}$ the reflector to its Birkhoff subcategory of abelian groups, and $\G$ the underlying set/free group comonad, it is well known that $H_{n} (A,\ab)_{\G}$ is just the $n$-th integral homology group $H_{n} (A,\Z)$ of $A$ (see page II.6.16 of~\cite{Quillen}).
\end{example}

\subsection{Comonads derived from a suitable adjunction}\label{Subsection-Regular-Comonads}
 From now on we shall assume that $\Ac$ is semi-abelian. We are going to construct simplicial resolutions for higher extensions: to do so, we consider a comonad $\G$ on $\Ac$ which is derived from an adjunction
\begin{equation}\label{comonadadjunction}
\vcenter{\xymatrix{
\Ac \ar@<1 ex>[r]^U \ar@{}[r]|{\top}  & \Xc, \ar@<1 ex>[l]^{F}}}
\end{equation}
satisfying conditions (\ref{Quillen1})--(\ref{Quillen2}) below. In this way, for every $n$, we shall obtain a comonad $\Gn$ on $\Arrn \Ac$ such that every $n$-extension $A$ in $\Ac$ yields a simplicial resolution $\Gn A$ \emph{in the category $\Extn \Ac$ of $n$-extensions}.
 
Let then $U \colon {\Ac \to \Xc}$ be a functor with a left adjoint $F \colon {\Xc \to \Ac}$ and write $\epsilon \colon {F \comp U \To 1_{\Ac}}$ and $\zeta \colon {1_{\Xc}\To U \comp F}$ for the counit and unit. This induces a comonad $\G$ as follows: $G=F \comp U$, $\epsilon$ is the counit and $\delta$ is the natural transformation defined by $\delta_{A}=F \zeta_{UA}$, for every object $A$ of $\Ac$. Let us now assume that the following conditions (given by Quillen on page II.5.5 of~\cite{Quillen}) are satisfied: 
\begin{enumerate}
\item\label{Quillen1}
$\epsilon_{A}\colon {GA\to A}$ is a regular epimorphism for all objects $A$ of~$\Ac$, and 
\item\label{Quillen2}
$FX$ is (regular epi)-projective for all objects $X$ of~$\Xc$. 
\end{enumerate} 
Then, in particular, $\epsilon_A\colon GA\to A$ is a \emph{projective presentation} of $A$, for every $A$ in~$\Ac$. 

\begin{example}\label{Example-Monadic}
Let $\Ac$ be semi-abelian and monadic over $\Xc = \Set$, and $\G$ the induced comonad on $\Ac$. Condition (\ref{Quillen1}) then follows from Beck's Theorem: actually, $\epsilon_{A}$ is a coequalizer of $G\epsilon_{A}$ and $\epsilon_{GA}$. And for every set $X$, $FX$ is free, hence (regular epi)-projective. 

Any variety of algebras is monadic over $\Set$, and thus semi-abelian varieties form an important class of examples. A characterization of such varieties is given by Bourn and Janelidze in their paper~\cite{Bourn-Janelidze}. In~\cite{Gran-Rosicky:Monadic}, Gran and Rosick\'y characterize semi-abelian categories, monadic over $\Set$.
\end{example}

\begin{definition}\label{Definition-Projective-extension}
An object of $\Arrn \Ac$ is called \emph{extension-projective} when it is projective with respect to the class $\Ec^{n}$ of $(n+1)$-extensions. An $n$-extension $f\colon {B\to A}$ is called a \emph{presentation} of $A$ when the object $B$ is extension-projective. An $n$-extension $f\colon {B\to A}$ is called an \emph{$n$-fold presentation} (or simply \emph{$n$-presentation}) when the object $B$ is extension-projective and $A$ is an $(n-1)$-presentation. (A \emph{$0$-presentation} is an object of $\Ac$.)
\end{definition}

\begin{lemma}
An object $f\colon {B\to A}$ of $\Arrn\Ac$ is extension-projective if and only if both $B$ and $A$ are extension-projective.\noproof
\end{lemma}

This lemma implies that an object $f\colon {B\to A}$ of $\Arrn\Ac$ is extension-projective if and only if all objects in the corresponding $n$-cube are projective in~$\Ac$.

\begin{proposition}\label{projectiveadjunctiongoesup}
Consider an adjunction \ref{comonadadjunction} that satisfies conditions (\ref{Quillen1})--(\ref{Quillen2}). Then, for every $n\geq 0$, there is an induced adjunction
\begin{equation}\label{comonadadjunctionn}
\xymatrix{
\Arrn\Ac \ar@<1 ex>[r]^-{U_n} \ar@{}[r]|{\top}  & \Arrn\Xc \ar@<1 ex>[l]^-{F_n}}
\end{equation}
defined degreewise, by putting, for objects $f\colon B\to A$ of $\Arrn\Ac$ and $z\colon Y\to X$ of $\Arrn\Xc$,
\begin{align*}
U_n(f\colon B\to A) &= U_{n-1}f\colon U_{n-1}B\to U_{n-1}A;\\
F_n(z\colon Y\to X) &= F_{n-1}z\colon F_{n-1}Y\to F_{n-1}X;\\
\epsilon_n (f\colon B\to A) &= (\epsilon_{n-1} B,\epsilon_{n-1} A)\colon F_nU_n f  \to f;\\
\zeta_{n} (z\colon Y\to X) &= (\zeta_{n-1} Y,\zeta_{n-1} X)\colon z \to U_n F_n z.
\end{align*}
We will often omit the indices and simply write $U$, $F$, $\epsilon$ and $\zeta$. Putting $G_n=F_n\comp U_n$ and $(\delta_n)_{A}=F( \zeta_n)_{U_nA}$ for every object $A$ of $\Arrn\Ac$ defines a comonad $\G_n=(G_n,\delta_n,\epsilon_n)$, denoted by $(G,\delta,\epsilon)$ when confusion is unlikely. Moreover, the following two conditions are satisfied:
\begin{enumerate}
\item
$\epsilon_{A}\colon {GA\to A}$ is an $(n+1)$-extension for all $n$-extensions $A$, and 
\item
$FX$ is extension-projective for all objects $X$ of~$\Arrn\Xc$. 
\end{enumerate} 
\end{proposition}
\begin{proof}
We use induction on $n$.

Since an object $f\colon {B\to A}$ of $\Arrn\Ac$ is extension-projective if and only if both $A$ and $B$ are extension-projective, condition (\ref{Quillen2}) follows immediately from the induction hypothesis. 

In order to show that condition (\ref{Quillen1}) holds as well, it suffices to prove that, for any $n$-extension $f\colon {B\to A}$, also $Gf\colon {GB\to GA}$ is an $n$-extension. Indeed, if this is the case, then the kernel pair $R[Gf]$ of the arrow $Gf$ is an $(n-1)$-extension, by Corollary \ref{Corollary-Higher-Denormalized-Rotlemma}. Also, $R[f]$ is an $(n-1)$-extension, being the kernel pair of the $n$-extension $f$. Hence, by the induction hypothesis, $\epsilon_{R[f]}$ is an $n$-extension. Furthermore, by the universal property of the kernel pair $R[Gf]$, there exists an arrow $a\colon {GR[f]\to R[Gf]}$ such that $\overline{\epsilon}_B \comp a=\epsilon_{R[f]}$, where $\overline{\epsilon}_B$ is the restriction of $\epsilon_B$ to ${R[Gf]\to R[f]}$. By Lemma \ref{Lemma-Higher-Ext-Compo} this implies that $\overline{\epsilon}_B$ is an $n$-extension. Finally, by Corollary \ref{Corollary-Higher-Denormalized-Rotlemma}, this implies that the right hand square in the next diagram is an $(n+1)$-extension.
\[
\xymatrix{
GR[f] \ar[d]_-{a} \ar@<1.2 ex>[rd] \ar@<0.2 ex>[rd] &&\\
R[Gf] \ar@<0.5 ex>[r] \ar@<-0.5 ex>[r] \ar[d]_-{\overline{\epsilon}_B} & GB \ar[r]^-{Gf} \ar[d]_-{\epsilon_B} & GA \ar[d]^-{\epsilon_A}\\
R[f] \ar@<0.5 ex>[r] \ar@<-0.5 ex>[r] & B \ar[r]_-f & A}
\]
Let us then prove that $Gf\colon {GB\to GA}$ is indeed an $n$-extension, for any $n$-extension $f\colon {B\to A}$. 

Let $X$ be any object of $\Arr^{n-1}\Xc$. Since, by assumption,  $FX$ is extension-projective, the map of $\hom$-sets
\[
\hom (FX, f) \colon \hom (FX, B) \to \hom(FX, A)
\]
is surjective. Since $F$ is left adjoint to $U$, this implies that the map
\[
\hom (X, Uf) \colon \hom(X,UB) \to \hom(X,UA)
\]
is surjective as well. Since this is the case for any object $X$ of $\Arr^{n-1}\Xc$, it follows that $Uf\colon UB\to UA$ is a split epimorphism. Consequently, also $Gf\colon GB\to GA$ is a split epimorphism. Since both its domain and codomain are $(n-1)$-extensions, it is an $n$-extension (see Example \ref{Example-Split-Epi-extension}).
\end{proof}

In particular, $\Extn \Ac$ has enough $\Ec^n$-projectives, because for any $n$-extension $A$, the map $\epsilon_A\colon GA\to A$ is an $(n+1)$-extension with an $\Ec^n$-projective domain.

\begin{corollary}\label{Corollary-Canonical-n-Fold-extension}\index{epsilonn@$\epsilon^{n}$}
If $\G$ is a comonad on a semi-abelian category $\Ac$ induced by an adjunction \ref{comonadadjunction} satisfying conditions (\ref{Quillen1})--(\ref{Quillen2}), then every object $A$ of $\Ac$ has a \emph{canonical $n$-presentation:} $\epsilon^{n}_{A}=\epsilon_{\epsilon_{A}^{n-1}}\colon {G\epsilon_{A}^{n-1}\to \epsilon_{A}^{n-1}}$.\noproof  
\end{corollary}

\section{Homology of extensions with respect to centralization}\label{Section-Homology-of-extensions}

 From now on, we concentrate on the situation where
\begin{enumerate}
\item $\Ac$ is a semi-abelian category and $\Bc$ is a Birkhoff subcategory of $\Ac$;
\item $\G$ is a comonad on $\Ac$ like in Subsection \ref{Subsection-Regular-Comonads}: $\G$ is derived from an adjunction \ref{comonadadjunction} which satisfies the conditions (\ref{Quillen1})--(\ref{Quillen2}) from Subsection \ref{Subsection-Regular-Comonads}.
\end{enumerate}
We showed that for all $n$, there is a Galois structure $\Gamma_n$ on $\Extn\Ac$ and an $\Ec^n$-comonad $\G_n$ on $\Arrn\Ac$. Let $f\colon {B\to A}$ be a presentation of an $(n-1)$-extension $A$ ($f$ may for instance be an $n$-presentation). In this section, we describe the homology with respect to $\Bc$-centralization of $f$ in terms of the homology with respect to $\CExt^{n-1}_\Bc\!\Ac$ of $A$. If $I$ denotes the reflector of $\Ac$ onto~$\Bc$, and $I_{1}\colon {\Ext \Ac\to \CExt_\Bc\Ac}$ the centralization functor, then for all $k\geq 2$, 
\[
H_{k} (f,I_{1})_{\G_{1}}\cong (H_{k+1} (A,I)_{\G}\to 0),
\] 
where $\G_{1}$ is the comonad on $\Arr \Ac$, induced by $\G$. Observe that the isomorphism above lives in the category $\Arr\Ac$.

More generally, if $I_{n}$ denotes the reflector of $\Extn\Ac$ onto $\CExt^{n}_\Bc\Ac$, then for all $n\geq 1$ and $k\geq 2$,
\[
H_{k} (f,I_{n})_{\Gn}\cong (H_{k+1} (A,I_{n-1})_{\G_{n-1}}\to 0),
\]
where $\Gn$ is the comonad on $\Arrn \Ac$, induced by $\G$.

\begin{remark}\label{Remark-Functor-I}
For any $n$, we shall consider $I_{n}$ as a functor ${\Extn \Ac \to \Arrn \Ac}$; then the situation fits into the scope of Definition~\ref{Definition-Homology}.
\end{remark}

\begin{lemma}\label{Lemma-I-of-I_{1}}
For every $n\geq 1$, given an $n$-extension $f$, there exists an isomorphism $I_{n-1}I_{n}f\cong I_{n-1}f$ in $\Arrn \Ac$.
\end{lemma}
\begin{proof}
Suppose $f\colon {B\to A}$ is an $n$-extension. Consider the $3\times 3$-diagram below. Note that all rows are short exact sequences, as is the middle column. 
\[
\xymatrix{
J_{n}[f] \ar@{=}[r] \ar[d] \ar@{}[rd]|{\texttt{(iv)}} & J_n[f] \ar@{{ >}->}[d] \ar@{-{>>}}[r] & 0\ar[d] \\
J_{n-1}B \ar@{{ >}->}[r] \ar[d]_-{J_{n-1}\rho^n_f} & B  \ar@{}[rd]|-{\texttt{(v)}} \ar@{-{>>}}[r] \ar@{-{>>}}[d]_-{\rho^n_f} & I_{n-1}B \ar[d] \\
J_{n-1}I_n[f] \ar@{{ >}->}[r]_-{\mu^{n-1}_{I_n[f]}} & I_n[f] \ar@{-{>>}}[r] & I_{n-1}I_n[f] }
\]
We must prove that the right hand column is short exact. By the strongly $\Ec^{n-1}$-Birkhoff property of $I_{n-1}$, the square \texttt{(v)} is an $n$-extension, hence $J_{n-1}\rho^n_f$ is an $(n-1)$-extension. In particular, it is a regular epimorphism. Since, furthermore, the square \texttt{(iv)} is a pullback (because $J_n[f]\subseteq J_{n-1}B$), the left hand column is short exact, hence, by applying the $3\times 3$-Lemma in the category ${\mathsf{Arr}^{n-1}\!} \Ac$, so is the right hand one.    
\end{proof}

\begin{lemma}\label{Lemma-I-vs-I_1}
Let $\G$ be a comonad on $\Ac$ like in Subsection \ref{Subsection-Regular-Comonads}. For every $n\geq 1$, given an $n$-extension $f\colon B\to A$, there exists an isomorphism $K[I_n\G_n f]\cong K[I_{n-1}\G_n f]$ in $\Sc\Arr^{n-1}\!\Ac$.
\end{lemma}
\begin{proof}
$I_{n}\Gn f$, which consists degreewise of central extensions with an $\Ec^n$-projective codomain, is degreewise trivial by Proposition~\ref{Central-is-normal}. This means that its kernel $K[I_{n}\Gn f]$ is isomorphic to $K[I_{n-1}I_{n}\Gn f]$:
\[
\resizebox{\textwidth}{!}
{\xymatrix{0 \ar[r] & K[I_{n}\Gn f] \ar@{.>}[d]_-{\cong} \ar@{{ >}->}[rr] && I_{n}[\Gn f] \ar@{-{>>}}[rr]^-{I_{n}\Gn f} \ar@{-{>>}}[d]_-{\eta^{n-1}_{I_{n}[\Gn f]}} \ar@{}[rrd]|<{\pullback} && \G_{n-1} A \ar@{-{>>}}[d]^-{\eta^{n-1}_{\G_{n-1} A}} \ar[r] & 0\\
0 \ar[r] & K[I_{n-1}I_{n}\Gn f] \ar@{{ >}->}[rr] && I_{n-1}I_{n}[\Gn f] \ar@{-{>>}}[rr]_-{I_{n-1}I_{n}\Gn f} && I_{n-1}\G_{n-1} A \ar[r] & 0}}
\]
By Lemma~\ref{Lemma-I-of-I_{1}}, the bottom short exact sequence is isomorphic to
\begin{equation}\label{Sequence-I_n}
%\resizebox{\textwidth}{!}
{\xymatrix{0 \ar[r] & K[I_{n-1}\G_{n} f] \ar@{{ >}->}[rr] && I_{n-1}\G_{n-1}B \ar@{-{>>}}[rr]^-{I_{n-1}\G_{n} f} && I_{n-1}\G_{n-1} A \ar[r] & 0,}}
\end{equation}
which implies that $K[I_n\G_n f]\cong K[I_{n-1}\G_n f]$.
\end{proof}

\begin{theorem}\label{Theorem-Homology-of-extension}
For any $n\geq 1$, $k\geq 2$ and any presentation $f\colon {B\to A}$ of an $(n-1)$-extension $A$,
\[
H_{k} (f,I_{n})_{\Gn}\cong (H_{k+1} (A,I_{n-1})_{\G_{n-1}}\to 0).
\]
\end{theorem}
\begin{proof}
We need to compute the homology of 
\[
I_{n}\Gn f\colon {I_{n}[\Gn f]\to \G_{n-1} A},
\]
a simplicial object in $\CExt^{n}_\Bc\Ac$---the image through $I_{n}$ of the simplicial resolution $\Gn f\colon {\G_{n-1} B\to \G_{n-1} A}$ in $\Extn \Ac$ of the extension $f$ of $A$. By definition, the homology objects $H_{k-1}I_{n}\Gn f$ may be computed degreewise, i.e.,
\[
H_{k-1}I_{n}\Gn f \cong (H_{k-1}I_{n}[\Gn f]\to H_{k-1}\G_{n-1}A).
\]
Note that $H_{i-1}\G_{n-1} A=0$ for all $i\geq 2$; by the long exact homology sequence~\cite[Corollary~5.7]{EverVdL2}, this implies that $H_{k-1}I_{n}[\Gn f]$ is isomorphic to $H_{k-1}K[I_{n}\Gn f]$. By Lemma~\ref{Lemma-I-vs-I_1}, we know that the simplicial objects $K[I_{n}\Gn f]$ and $K[I_{n-1}\G_n f]$ are isomorphic. Using the long exact homology sequence associated with Sequence~\ref{Sequence-I_n} and the fact that $B$ is extension-projective (so that ${\G_{n-1}B\to B}$ is contractible), the desired isomorphism 
\[
H_{k-1}I_{n}[\Gn f]\cong H_{k}I_{n-1}\G_{n-1} A
\]
is obtained.
\end{proof}

\section{The Hopf formula for the second homology object}\label{Section-Second-Homology-Object}

In this section we give a direct proof of the Hopf formula from~\cite{EverVdL2}, which describes the second homology object in terms of ``generalized commutators".

Recall~\cite[Corollary~3.10]{EverVdL2} that $H_{0}S=\Coeq [\del_{0},\del_{1}\colon {S_{1}\to S_{0}}]$ for every simplicial object $S$ in a semi-abelian category $\Ac$---a consequence of the fact that in $\Ac$, a regular epimorphism is a cokernel of its kernel. 

\begin{lemma}\label{Lemma-H_0}
Let $f\colon B\to A$ be an $n$-extension; then $H_0I_n\G_n f\cong I_n f$.
\end{lemma}
\begin{proof}
The diagram
\[
\xymatrix{I_n[G^2f] \ar@<.5ex>[r]^-{I_n[G\epsilon_f]} \ar@<-.5ex>[r]_-{I_n[\epsilon_{Gf}]} \ar@{-{ >>}}[d]_-{I_nG^2f} & I_n[Gf] \ar@{-{ >>}}[r]^-{I_n[\epsilon_{f}]} \ar@{-{ >>}}[d]^-{I_nGf} & I_n[f] \ar@{-{ >>}}[d]^-{I_nf} \\
G^2A \ar@<.5ex>[r]^-{G\epsilon_A} \ar@<-.5ex>[r]_-{\epsilon_{GA}} & GA \ar@{-{ >>}}[r]_-{\epsilon_A} & A}
\]
is a coequalizer, because the square 
\[
\xymatrix{GB \ar@{-{ >>}}[d]_-{\rho^n_{Gf}} \ar@{-{ >>}}[r]^-{\epsilon_B} & B \ar@{-{ >>}}[d]^-{\rho^n_f} \\
I_n[Gf] \ar@{-{ >>}}[r]_-{I_n[\epsilon_{f}]} & I_n[f]}
\]
is a pushout: this follows from the strongly $\Ec^n$-Birkhoff property of $I_n$.
\end{proof}

\begin{theorem}\label{Theorem-Hopf-2}\cite{EverVdL2} Let $f\colon {B\to A}$ be a $1$-presentation; then 
\[
H_{2} (A,I)_{\G}\cong\frac{JB\cap K[f]}{J_{1}[f]}.
\]
\end{theorem}
\begin{proof}
Note that $J_{1}[f]$ is a normal subobject of $JB\cap K[f]$ by definition. Since $B$ is extension-projective, $I\G B$ is contractible, hence the short exact sequence of simplicial objects
\[
\xymatrix{0 \ar[r] & K[I\G f] \ar@{{ >}->}[rr]^-{\ker I\G f} && I\G B \ar@{-{>>}}[rr]^-{I\G f} && I\G A \ar[r] & 0}
\]
yields an exact homology sequence
\[
\xymatrix{
0 \ar[r] & H_2(A,I)_{\G}=H_1I\G A \ar@{{ >}->}[r] & H_0K[I\G f] \ar[r] & IB \ar@{-{>>}}[r]^-{If} & IA \ar[r] & 0.}
\]
Taking into account Lemma~\ref{Lemma-I-vs-I_1}, we thus get a short exact sequence
\begin{equation}\label{shortexacthopf}
\xymatrix{
0 \ar[r] & H_2(A,I)_{\G} \ar@{{ >}->}[r] & H_0K[I_1\G f] \ar@{-{>>}}[r] & K[If] \ar[r] & 0.}
\end{equation}
Now, using the long exact homology sequence induced by
\[
\xymatrix{0 \ar[r] & K[I_{1}\G f] \ar@{{ >}->}[rr] && I_{1}[\G f] \ar@{-{>>}}[rr]^-{I_{1}\G f} && \G A \ar[r] & 0}
\]
we find that 
\[
H_{0}K[I_{1}\G f] = K[H_0I_1 \G f]. 
\]
Hence, by the foregoing lemma,
\[
H_{0}K[I_{1}\G f] =  K[I_{1}f].
\]
By the $3\times 3$-Lemma,
\[
K\Bigl[I_{1}f\colon  \frac{B}{J_1[f]}\to A\Bigr] = \frac{K[f]}{J_1[f]}
\]
and
\[
K\Bigl[If\colon \frac{B}{JB} \to \frac{A}{ JA}\Bigr] = \frac{K[f]}{JB\cap K[f]}.
\]
The theorem now follows from applying Noether's First Isomorphism Theorem (Theorem~4.3.10 in~\cite{Borceux-Bourn}, a direct consequence of the $3 \times 3$-Lemma) to~\ref{shortexacthopf}.
\end{proof}

This also works for centralization of $n$-extensions:

\begin{proposition}\label{Proposition-Hopf-2}
Consider $n\in \N$ and $f\colon {B\to A}$, a presentation of an $n$-extension~$A$. Then
\[
H_{2} (A,I_{n})_{\Gn}\cong \frac{J_nB\cap K[f]}{J_{n+1}[f]}:
\]
$H_{2} (A,I_{n})_{\Gn}$ is the direct image of $K[f]\cap J_nB$ along $\rho^{n+1}_{f}\colon {B\to I_{n+1}[f]}$.\noproof 
\end{proposition}

\begin{remark}\label{Remark-Baer-Invariants}
In particular, the expressions on the right hand side of these Hopf formulae are \emph{Baer invariants}~\cite{Froehlich, EverVdL1}: they are independent of the chosen presentation of $A$.
\end{remark}

\begin{remark}\label{Remark-Dimension-in-Hopf-2}
Note that the objects in this Hopf formula are in $\Arrn \Ac$.
\end{remark}

\begin{notation}\label{Notation-Higher-Extension}\index{fA@$f_{A}$}
The functor category $\Arrn \Ac =\Hom (\Two^{n},\Ac)$ may with advantage be described in the following way~\cite{Brown-Ellis}. Let $\set{n}$ denote the set $\{1,\dots,n  \}$, $\set{0}=\varnothing$. Then the category $\Two^{n}$ is isomorphic to the power-set $\Pc \set{n}$, the set of subsets of $\set{n}$, ordered by inclusion. This means that an inclusion $A\subseteq B$ in $\set{n}$ corresponds to a map $\iota^{A}_{B}\colon {A\to B}$ in $\Two^{n}$. We write $f_A$ for $f(\iota^{A}_{A})$ and $f^A_B$ for $f(\iota^A_B)$. When, in particular, $A$ is $\varnothing$ and $B$ is a singleton $\{i \}$, we write $f_{i}=f^{\varnothing}_{\{i \}}$. For instance, $f_{\varnothing}$ is the ``initial object'' of the cube
\[ 
\xymatrix@!=.5cm{& f_{\varnothing} \ar@{-{ >>}}[dd]|-{\hole}_(.75){f_{1}} \ar@{-{ >>}}[rr]^-{f_{3}} \ar@{-{ >>}}[ld]_-{f_{2}} && f_{\{3\}} \ar@{-{ >>}}[dd]^-{f^{\{3\}}_{\{1,3\}}} \ar@{-{ >>}}[ld]_-{f^{\{3\}}_{\{2,3\}}}\\
f_{\{2\}}  \ar@{-{ >>}}[dd]_-{f^{\{2\}}_{\set{2}}} \ar@{-{ >>}}[rr] && f_{\{2,3\}} \ar@{-{ >>}}[dd]\\
& f_{\{1\}} \ar@{-{ >>}}[rr]|(.5){\hole} \ar@{-{ >>}}[ld] && f_{\{1,3\}} \ar@{-{ >>}}[ld]\\
f_{\{1,2\}} \ar@{-{ >>}}[rr]_-{f^{\{1,2\}}_{\set{3}}} && f_{\{1,2,3\}}}
\]
and $f_{i}$, for $i\in \set{3}$, is any of its three initial ribs. 

Let us now write the Hopf formula using the notations just defined. By Notation~\ref{Notation-W}, the kernel $K[\eta^{n}_{B}]= J_n B$ is written $\Psi^{n}L_{n}[B]$; it follows that the intersection $J_n B \cap K[f]$ is $\Psi^{n}(L_{n}[B]\cap K[f_{n+1}])$. Also, the kernel  $J_{n+1}[{f}]$ is $\Psi^{n}L_{n+1}[f]$, and hence
\[
H_{2} (A,I_{n})_{\Gn}\cong \Psi^{n}\frac{L_{n}[B]\cap K[f_{n+1}]}{L_{n+1}[f]}.
\]
\end{notation}

\begin{notation}\label{Notation-Lower-from-Higher}\index{extnf@$\ext_{n}f$}
An $n$-extension $f$ naturally gives rise to $k$-extensions for all $0\leq k\leq n$. We make the following choice of $\ext_{k}f\in \Ext^{k}\!\Ac$: put $\ext_{n}f=f$ and $\ext_{k-1}f=\cod\ext_{k}f$ ($0< k\leq n$). It is easily seen that $\ext_{k}f$ may be obtained from $f$ by precomposing with the functor $\Two^{k}\to \Two^{n}$ determined by $A\mapsto A\cup (\set{n}\setminus\set{k})$. 
\end{notation}

\section{The higher Hopf formulae}\label{Section-Higher-Hopf-Formulae}

Let us start this section on the higher Hopf formulae by explaining the formula for the third homology object in the special case of canonically chosen presentations; the others---where $n>3$ or presentations are chosen non-canonically---work by essentially the same principle. 

Suppose that $\Ac$ is a semi-abelian category and $\Bc$ is Birkhoff in $\Ac$. Let $\G$ be a comonad on $\Ac$ derived from an adjunction \ref{comonadadjunction} which satisfies the conditions (\ref{Quillen1})--(\ref{Quillen2}) from Subsection \ref{Subsection-Regular-Comonads}. Consider an object $A$ in $\Ac$. By Theorem~\ref{Theorem-Homology-of-extension},
\[
\bigl(H_{3} (A,I)_{\G}\to 0\bigr) \cong H_{2} (\epsilon_{A},I_{1})_{\G_{1}},
\] 
i.e., in order to compute the third homology object of $A$ with respect to $I$, it suffices to compute the second homology extension with respect to $I_{1}$ of the canonical presentation $\epsilon_{A}\colon {GA\to A}$. This may be done using Proposition~\ref{Proposition-Hopf-2}:
\[
H_{2} (\epsilon_{A},I_{1})_{\G_{1}} \cong \rho^{2}_{\epsilon_{A}^{2}} (K[\epsilon_{A}^{2}]\cap K[\eta^{1}_{G\epsilon_{A}}]).
\]
Here $\epsilon_{A}^{2}\colon G\epsilon_A\to \epsilon_A$ is the canonical double extension of $A$, which has as kernel $K[\epsilon^{2}_{A}]$:
\[
\xymatrix{0 \ar[r] & K[\epsilon_{GA}] \ar@{-{ >>}}[d]_-{K[\epsilon_{A}^{2}]} \ar@{{ |>}->}[r] & G^{2}A  \ar@{-{ >>}}[d]_-{G\epsilon_{A}} \ar@{-{ >>}}[r]^-{\epsilon_{GA}} & GA \ar@{-{ >>}}[d]^-{\epsilon_{A}} \ar[r] & 0\\
0 \ar[r] & K[\epsilon_{A}] \ar@{{ |>}->}[r] & GA \ar@{-{ >>}}[r]_-{\epsilon_{A}} & A \ar[r] & 0.}
\]
The other factor in the intersection is the kernel $K[\eta^{1}_{G\epsilon_{A}}]$ of the right hand side square in the next diagram.
\[
\xymatrix{0 \ar[r] & J_{1}[G\epsilon_{A}] \ar@{-{ >>}}[d]_-{J_{1}G\epsilon_{A}=K[\eta^{1}_{G\epsilon_{A}}]} \ar@{{ |>}->}[r] & G^{2}A  \ar@{-{ >>}}[d]_-{G\epsilon_{A}} \ar@{-{ >>}}[r] & I_{1}[G\epsilon_{A}] \ar@{-{ >>}}[d]^-{I_{1}G\epsilon_{A}} \ar[r] & 0\\
0  \ar[r] & 0 \ar@{{ |>}->}[r] & GA \ar@{=}[r] & GA \ar[r] & 0}
\]
The map $\rho^{2}_{\epsilon_{A}^{2}}\colon {G\epsilon_{A}\to I_2\epsilon_{A}^{2}}$ may be pictured as
\[
\xymatrix{G^{2}A \ar@{-{ >>}}[d]_{G\epsilon_{A}} \ar@{-{ >>}}[r]^-{\eta_{\epsilon_{A}^{2}}} & I_{2}[\epsilon_{A}^{2}]_{0} \ar@{-{ >>}}[d]\\
GA \ar@{=}[r] & GA.}
\]
Hence the direct image of $K[\eta^{1}_{G\epsilon_{A}}]\cap K[\epsilon_{A}^{2}]$ along $\rho^{2}_{\epsilon_{A}^{2}}$ is the front square of the cube 
\[
\xymatrix@!=.5cm{&  J_{1}[G\epsilon_{A}]\cap K[\epsilon_{GA}] \ar@{-{ >>}}[dd]|(.5){\hole} \ar@{{ |>}->}[rr] \ar@{-{ >>}}[ld] && G^{2}A \ar@{-{ >>}}[dd]^-{G\epsilon_{A}} \ar@{-{ >>}}[ld]\\
\frac{J_{1}[G\epsilon_{A}]\cap K[\epsilon_{GA}]}{K[\rho^{2}_{\epsilon_{A}^{2}}]} \ar@{-{ >>}}[dd] \ar@{{ |>}->}[rr] && I_{2}[\epsilon_{A}^{2}]_{0} \ar@{-{ >>}}[dd]\\
& 0 \ar@{{ |>}->}[rr]|(.5){\hole} \ar@{=}[ld] && {GA} \ar@{=}[ld]\\
0 \ar@{{ |>}->}[rr] && GA.}
\]
If we write $L_{2}[\epsilon_{A}^{2}]=K[\rho^{2}_{\epsilon_{A}^{2}}]$, thus we get that 
\[
H_{3} (A,I)_{\G} \cong \frac{J_{1}[G\epsilon_{A}]\cap K[\epsilon_{GA}]}{L_{2}[\epsilon_{A}^{2}]}.
\]
Using the Hopf formula for the second homology object, this may be rewritten as
\begin{equation}\label{Latter-form}
H_{3} (A,I)_{\G} \cong \frac{JG^{2}A\cap K[G\epsilon_{A}]\cap K[\epsilon_{GA}]}{L_{2}[\epsilon_{A}^{2}]},
\end{equation}
because $GA$ is projective, and hence $H_{2} (GA,I)_{\G}=0$, which implies that
\[
J_{1}[G\epsilon_{A}]\cong JG^{2}A\cap K[G\epsilon_{A}].
\]
Modulo the more abstract denominator, this latter form~\ref{Latter-form} of the formula is how it occurs in Brown and Ellis's paper~\cite{Brown-Ellis}. The main difference is that here, the formula is also valid for other categories $\Ac$ than the category $\Gp$ of groups, and other Birkhoff subcategories $\Bc$ than the category $\Ab$ of abelian groups.

Following the lines of the method due to Janelidze as set out in~\cite{Janelidze:Double} and sketched above, we now show that also the higher Hopf formulae from~\cite{Brown-Ellis} may be generalized to arbitrary semi-abelian categories. 

We speak of {\it an $n$-presentation $f$ of an object $A$} in $\Ac$ if $f$ is an $n$-presentation and if $f_{\set{n}} = A$ or, in other words, if $\ext_{0}f=A$.

\begin{theorem}\label{Theorem-Hopf-BI}
Let $f$ be an $n$-presentation of an object $A$ of $\Ac$. Then
\begin{equation}\label{Formula-Hopf-BI}
H_{n+1} (A,I)_{\G}\cong \frac{Jf_{\varnothing}\cap \bigcap_{i\in \set{n}}K[f_{i}]}{L_{n}[f]}.
\end{equation}
In particular, the expression on the right hand side of Formula~\ref{Formula-Hopf-BI} is a (higher) Baer invariant: independent of the chosen $n$-presentation of $A$. 
\end{theorem}  
\begin{proof}
The formula holds by induction on $n$. Consider the following chain of isomorphisms (cf.\ Notation~\ref{Notation-Lower-from-Higher}):
\begin{align*}
\Psi^{n-1}H_{n+1} (A,I)_{\G} &\cong \Psi^{n-2}H_{n} (\ext_{1}f,I_{1})_{\G_{1}}\\
 &\cong \cdots \cong \Psi H_{3} (\ext_{n-2}f,I_{n-2})_{\G_{n-2}}\\
 &\cong H_{2} (\ext_{n-1}f,I_{n-1})_{\G_{n-1}}\\
 &\cong \Psi^{n-1}\frac{ L_{n-1}[\dom f] \cap K[f_{n}]}{L_{n}[f]}\\
 &\cong \Psi^{n-1}\frac{Jf_{\varnothing}\cap \bigcap_{i\in \set{n}}K[f_{i}]}{L_{n}[f]}.
\end{align*}
The first isomorphism follows from Theorem~\ref{Theorem-Homology-of-extension}. Repeatedly applying this theorem gives the expression containing $H_{2}$. Proposition~\ref{Proposition-Hopf-2} and Notation~\ref{Notation-Higher-Extension} account for the next isomorphism, and the last one follows by the induction hypothesis, the Hopf formula for $H_{n} (f_{\set{n}\setminus\{n \}},I)_{\G}$, because $\dom f$ is an $(n-1)$-presentation of $f_{\set{n}\setminus\{n \}}$. Indeed, $f_{\set{n}\setminus\{n \}}$ is a projective object, and hence $H_{n} (f_{\set{n}\setminus\{n \}},I)_{\G}=0$, which implies that an isomorphism
\[
L_{n-1}[\dom f]\cong Jf_{\varnothing}\cap \bigcap_{i\in \set{n-1}}K[f_{i}]
\] 
exists.
\end{proof}
\begin{remark}\label{RemarkInd}
In particular, the formula in the Theorem above shows that the homology objects are independent of the chosen comonad on $\Ac$.
\end{remark}
When a canonical presentation of $A$ is chosen, one obtains the following formula.

\begin{corollary}\label{Corollary-Hopf}
For every $n\geq 1$, an isomorphism
\[
H_{n+1} (A,I)_{\G}\cong \frac{JG^{n}A\cap \bigcap_{i\in [n-1]}K[G^{n-1-i}\epsilon_{G^{i}A}]}{L_{n}[\epsilon_{A}^{n}]}
\]
exists.\noproof 
\end{corollary}

\section{Examples}
In this final section we consider some specific cases of Theorem~\ref{Theorem-Hopf-BI} and give explicit descriptions of the right hand side of the Hopf formulae. In particular, we will explain why Theorem~\ref{Theorem-Hopf-BI} gives Brown and Ellis's formulae~\cite{Brown-Ellis} in the case where $\Ac$ is $\Gp$, the variety of groups and $\Bc$ is $\Ab$, the subvariety of all abelian groups. More generally, if $\Ac$ is $\Gp$ and $\Bc$ is $\Nil_k$, the subvariety of $k$-nilpotent groups (for some positive integer~$k$), then we obtain the formulae due to Donadze, Inassaridze and Porter~\cite{Donadze-Inassaridze-Porter}. Furthermore, we find new formulae when $\Bc=\Sol_k$, the subvariety of $k$-solvable groups. Also, we explain that, by using the same arguments, similar formulae can be obtained from Theorem~\ref{Theorem-Hopf-BI} when $\Ac$ is, e.g., the variety of rings or the variety of Lie algebras. Finally, we shall consider the case where $\Ac$ is $\PXM$, the variety of precrossed modules, and $\Bc$ is $\XM$, the subvariety of crossed modules.

From now on, we shall drop $\G$ in the notation for homology. This is justified by Remark~\ref{RemarkInd} and by the fact that we shall be dealing with varieties and their canonical comonad.

\subsection{Groups vs.\ abelian groups}\label{abelianization}
Let us consider the Hopf formulae~\ref{Theorem-Hopf-BI} in the particular case where $\Ac=\Gp$ is the variety of groups and $\Bc=\Ab$ the subvariety of abelian groups. In order to simplify the arguments below, we define the commutator of two (not necessarily normal) subgroups $A$ and $B$ of a group $G$ as the \emph{normal} subgroup of $G$ generated by all elements $[a,b]=aba^{-1}b^{-1}$, with $a\in A$ and $b\in B$. Note that the reflector $I=\ab\colon \Gp\to\Ab$ is given by $\ab(G)=G/[G,G]$. It is known (see, e.g.,~\cite{Furtado-Coelho} and Example~\ref{centralextensionsforgroupsxample}) that in this situation $L_1[f]=[K[f],B]$ for any extension $f\colon {B\to A}$, hence when $f$ is a presentation, Formula~\ref{Theorem-Hopf-BI} yields, for $n=1$, the classical Hopf formula
\[
H_2(G,\ab) \cong \frac{[B,B] \cap K[f]}{[K[f],B]}.
\]
We are now going to show that our formulae for $H_n(G,\ab)$ coincide with those of Brown and Ellis, also when $n\geq 2$. In particular, we shall prove that, for all $n\geq 0$ and any $n$-extension $f$, 
\begin{equation}\label{downtheclassicalhopf}
L_n[f]= \prod_{I \subseteq \set{n}}\Bigl[ \bigcap_{i\in I} K[f_i],   \bigcap_{i\notin I} K[f_i] \Bigr],
\end{equation}
where it is understood that $ \bigcap_{i\in \varnothing} K[f_i]=f_{\varnothing}$. 

Let us recall the following \emph{Witt-Hall identities}, valid in any group $G$, for any elements $a$, $a_1$, $a_2$, $b$, $b_1$, $b_2\in G$:
\begin{eqnarray}
[a,b_1 b_2] & = & [a,b_1] \   {b_1}[a,b_2]b_1^{-1} \label{WH1} \\
%[a_1 a_2, b] 
{[}a_1a_2,b]
& = & {a_1}[a_2,b]a_1^{-1} \ [a_1,b]. \label{WH2}
\end{eqnarray}
Let us write $A\cdot B$ for the product of subgroups $A$ and $B$ of a group $G$. The identities~\ref{WH1} and~\ref{WH2} imply, for all groups $G$ and subgroups $A$, $A_1$, $A_2$, $B$, $B_1$ and $B_2$ of $G$ that
\[
[A,B_1\cdot B_2] = [A,B_1] \cdot [A,B_2]
\]   
and
\[
[A_1\cdot A_2, B] = [A_1,B] \cdot [A_2,B].
\]
For a group $G$, we denote by $\Delta_G=(1_{G},1_{G})\colon {G\to G\times G}$ the \emph{diagonal} of $G$: the smallest internal equivalence relation on $G$.

Then, in particular, we get

\begin{lemma}\label{joingroup}
Suppose $G$ is a group, $M$ and $N$ are normal subgroups of $G$ and $B$ is a subgroup of $G\times G$. Then the following identity holds:
\[
[M\times_{\frac{M}{M\cap N}} M, B] = [0\times (M\cap N), B] \cdot [\Delta_M, B].
\]
\end{lemma}
\begin{proof}
This follows from the above, by taking into account that 
\[
M\times_{\frac{M}{M\cap N}} M = 0\times (M\cap N) \cdot \Delta_M,
\]
for any normal subgroups $M$ and $N$ of a group $G$.
\end{proof}

Let us then prove the identity~\ref{downtheclassicalhopf}, for all $n\geq 0$ and any $n$-extension $f$. We do this by induction. For $n=0$, we have
\[
L_0[A]=JA=[A,A]= \prod_{I \subseteq \varnothing}[ A,  A ],
\]
when $A=f$ is a $0$-extension in $\Gp$. Now, assume that~\ref{downtheclassicalhopf} holds for some $n-1$. To show that~\ref{downtheclassicalhopf} holds for $n$, it suffices to prove that 
\begin{eqnarray}\label{zeroatthesametime}
L_n[f]=0 &\Leftrightarrow &\prod_{I \subseteq \set{n}}\Bigl[ \bigcap_{i\in I} K[f_i],   \bigcap_{i\notin I} K[f_i] \Bigr]=0,
\end{eqnarray}
for every $n$-extension $f$. Indeed, suppose the equivalence~\ref{zeroatthesametime} holds. Recall that, for any group $G$ and any subgroups $A$ and $B$ of $G$, $[A,B]$ is the smallest normal subgroup $N$ of $G$ such that 
\[
\Bigl[\frac{A\cdot N}{N}, \frac{B\cdot N}{N} \Bigr] = 0.
\]
 From this one deduces that 
\[
\prod_{I \subseteq \set{n}}\Bigl[ \bigcap_{i\in I} K[f_i],   \bigcap_{i\notin I} K[f_i] \Bigr]
\]
is the smallest normal subgroup $N$ of $f_{\varnothing}$ satisfying 
\[
\prod_{I \subseteq \set{n}}\Bigl[  \frac{\bigcap_{i\in I}K[f_i]\cdot N}{N},  \frac{ \bigcap_{i\notin I} K[f_i] \cdot N}{N}\Bigr]=0.
\]
Since $L_n[I_nf]=0$, applying~\ref{zeroatthesametime} to the $n$-extension $I_nf$ gives that
\[
\prod_{I \subseteq \set{n}}\Bigl[ \bigcap_{i\in I} \frac{K[f_i]}{L_n[f]},   \bigcap_{i\notin I} \frac{K[f_i]}{L_n[f]} \Bigr] = \prod_{I \subseteq \set{n}}\Bigl[ \bigcap_{i\in I} K[(I_nf)_i],   \bigcap_{i\notin I} K[(I_nf)_i] \Bigr]= 0, 
\]
hence 
\[
L_n[f]  \supseteq   \prod_{I \subseteq \set{n}}\Bigl[ \bigcap_{i\in I} K[f_i],   \bigcap_{i\notin I} K[f_i] \Bigr].
\]
The other inclusion follows similarly, now taking into account that $I_n$ is a reflector.

Let us then prove the equivalence~\ref{zeroatthesametime}. Suppose $f\colon {B\to A}$ is an $n$-extension. It follows readily from the definition of $J_{n}=\Psi L_{n}$ that $L_n[f]=0$ if and only if
\[
\pi_1L_{n-1}[R[f]] = \pi_2L_{n-1}[R[f]],
\]
where $(\pi_1,\pi_2)$ denotes the kernel pair of $f_n\colon {f_{\varnothing}\to f_{\{n\}}}$.
By the induction hypothesis,
\begin{eqnarray*}
L_{n-1}[R[f]] &=& \prod_{I \subseteq \set{n-1}}\Bigl[ \bigcap_{i\in I} K[R[f]_i], \bigcap_{i\notin I} K[R[f]_i] \Bigr].
\end{eqnarray*}
Recall from Notation~\ref{Notation-Higher-Extension} that $R[f]_A=R[f^A_{A\cup \{n\}}]$ for all $A\subseteq \set{n-1}$. In particular,  $R[f]_{\varnothing}=R[f_n]$ and $R[f]_{\{i\}}=R[f^{\{i\}}_{\{i,n\}}]$. Furthermore, $R[f]_i$ is such that the following diagram commutes. 
\[
\xymatrix{
R[f]_{\varnothing}
 \ar@<-0.5 ex>[r] \ar@<0.5 ex>[r] \ar[d]_{R[f]_i} & f_{\varnothing} \ar[r]^{f_n} \ar[d]_{f_i} & f_{\{n\}} \ar[d]^{f^{\{n\}}_{\{i,n\}}} \\
R[f]_{\{i\}} \ar@<-0.5 ex>[r] \ar@<0.5 ex>[r] & f_{\{i\}} \ar[r]_{f^{\{i\}}_{\{i,n\} } } & f_{\{i,n\}} }
\]
Since the right hand square is a regular pushout, Proposition~\ref{Proposition-Rotlemma-Kernels} implies that
\[
K[f^{\{n\}}_{\{i,n\}}]= \frac{K[f_i]}{K[f_n]\cap K[f_i]}.
\]
Consequently,
\[
K[R[f]_i] = K[f_i]\times_{\frac{K[f_i]}{K[f_n]\cap K[f_i]}} K[f_i].
\]
It follows that 
\begin{equation*}\label{kernelL_n}
L_{n-1}[R[f]] 
 =  \prod_{I \subseteq \set{n-1}}  \Bigl[ \bigcap_{i\in I} K[f_i]\times_{\frac{K[f_i]}{K[f_i]\cap K[f_n]}} K[f_i]    , \bigcap_{i\notin I} K[f_i]\times_{\frac{K[f_i]}{K[f_i]\cap K[f_n]}} K[f_i]     \Bigr].
\end{equation*}
Moreover note that, for all $I\subseteq \set{n}$, we have
\[
\bigcap_{i\in I} K[f_i]\times_{\frac{K[f_i]}{K[f_i]\cap K[f_n]}} K[f_i]   =    \bigcap_{i\in I} K[f_i]\times_{\frac{ \bigcap_{i\in I} K[f_i]}{ \bigcap_{i\in I} K[f_i]\cap K[f_n]}}  \bigcap_{i\in I} K[f_i].   
\]
Now, applying Lemma~\ref{joingroup}, we find that
\begin{eqnarray*}
 L_{n-1}[R[f]] 
 &=& 
  \prod_{I \subseteq \set{n-1}}  \Bigl[ 0 \times  \bigcap_{i\in I} K[f_i] \cap K[f_n], 0  \times  \bigcap_{i\notin I} K[f_i] \cap K[f_n] \Bigr] \\
 && \cdot   \Bigl[ 0 \times  \bigcap_{i\in I} K[f_i] \cap K[f_n],  \Delta_{\bigcap_{i\notin I} K[f_i]  }\Bigr] \\
 &&  \cdot \Bigl[ \Delta_{\bigcap_{i\in I} K[f_i]  } , 0  \times  \bigcap_{i\notin I} K[f_i] \cap K[f_n] \Bigr]\\
 &&  \cdot \Bigl[ \Delta_{\bigcap_{i\in I} K[f_i]  } , \Delta_{\bigcap_{i\notin I} K[f_i]  }   \Bigr].
\end{eqnarray*}
Consequently, $\pi_1L_{n-1}[R[f]]=\pi_2L_{n-1}[R[f]] $ if and only if
\[
\Bigl[ \bigcap_{i\in I} K[f_i],   \bigcap_{i\notin I} K[f_i]\cap K[f_n] \Bigr] = \{1\},
\]
for all $I\subseteq \set{n-1}$, i.e., if and only if
\[
 \prod_{I \subseteq \set{n}}\Bigl[ \bigcap_{i\in I} K[f_i],   \bigcap_{i\notin I} K[f_i] \Bigr]=\{1\},
\] 
which is exactly what we wanted. We conclude that the equality~\ref{downtheclassicalhopf} holds for all $n\geq 0$. Thus Theorem~\ref{Theorem-Hopf-BI} induces Brown and Ellis's formulae:

\begin{theorem}\label{Brown-Ellis}\cite{Brown-Ellis}
For any $n\geq 1$ and any $n$-presentation $f$ of a group $G$, an isomorphism
\begin{equation*}
H_{n+1} (G,\ab) \cong \frac{[f_{\varnothing}, f_{\varnothing}]  \cap \bigcap_{i\in \set{n}}K[f_{i}]}{  \prod_{I \subseteq \set{n}}\Bigl[ \bigcap_{i\in I} K[f_i],   \bigcap_{i\notin I} K[f_i] \Bigr]  }
\end{equation*}
exists.\noproof 
\end{theorem}

\subsection{Groups vs.\ $k$-nilpotent groups}
Recall that a group $A$ is $k$-\emph{nilpotent} (=~nilpotent of class at most $k$) if and only if $Z_kA= \{ 1\}$, where $Z_kA$ is the $k$-th term in the \emph{descending central series} of $A$ defined by 
\begin{eqnarray*}
Z_1 A& = & [A,A] ,\\
Z_2 A & = & [[A,A], A] ,\\
\cdots & = & \cdots \\
Z_k A & = & [Z_{k-1} A, A].
\end{eqnarray*}
Let us now consider the Hopf formulae~\ref{Theorem-Hopf-BI} in the situation where $\Ac=\Gp$ and $\Bc=\Nil_k$, the subvariety of $\Gp$ of all $k$-nilpotent groups (for some positive integer~$k$). Let us write $\nil_k$ for the reflector $\Gp\to\Nil_k$, which sends a group $A$ to the quotient $A/Z_k A$.
 
The arguments used above in Subsection~\ref{abelianization} can easily be adapted to show that Lemma~\ref{joingroup}  implies, in this situation, that
\[
L_n[f]=    \prod_{I_1 \cup \dots \cup I_k=\set{n}} \Bigl[\dots \Bigl[\Bigl[ \bigcap_{i\in I_1} K[f_i],   \bigcap_{i\in I_2} K[f_i] \Bigr],   \bigcap_{i\in I_3} K[f_i] \Bigr], \dots \Bigr], 
\]
for any $n\geq 0$ and any $n$-extension of groups $f$. Consequently, the Hopf formulae~\ref{Theorem-Hopf-BI} yield
 
\begin{theorem}\label{Donadze}\cite{Donadze-Inassaridze-Porter}
For any $n\geq 1$ and any $n$-presentation $f$ of a group $G$, an isomorphism
\[
%\resizebox{\textwidth}{!}
{\displaystyle H_{n+1} (G,\nil_k) \cong 
\frac{[ \dots [[f_{\varnothing}, f_{\varnothing}], f_{\varnothing}],\dots ] \cap \bigcap_{i\in \set{n}}K[f_{i}]}{  \prod_{I_1 \cup \dots \cup I_k=\set{n}} \Bigl[\dots \Bigl[\Bigl[ \bigcap_{i\in I_1} K[f_i],   \bigcap_{i\in I_2} K[f_i] \Bigr],   \bigcap_{i\in I_3} K[f_i] \Bigr], \dots \Bigr] }}
\]
exists.\noproof 
\end{theorem}

\subsection{Groups vs.\ $k$-solvable groups}
Recall that a group $A$ is $k$-\emph{solvable} (=~solvable of class at most $k$) if and only if $D_k A= \{ 1\}$, where $D_k A$ is the $k$-th term in the \emph{derived series} of $A$ defined by 
\begin{eqnarray*}
D_1 A& = & [A,A], \\
D_2 A & = & [[A,A], [A,A]], \\
\cdots & = & \cdots \\
D_k A & = & [D_{k-1} A, D_{k-1}A].
\end{eqnarray*}
Let us consider the Hopf formulae~\ref{Theorem-Hopf-BI} in the case where $\Ac=\Sol_k$, the subvariety of $\Gp$ of all $k$-solvable groups (for some positive integer~$k$). We write $\sol_k$ for the reflector $\Gp\to\Sol_k$, which sends a group $A$ to the quotient $A/D_k A$. In order to express a formula for $L_n[\cdot]$, it is useful to introduce also the following notation:

Given normal subgroups $X_1$, $X_2$, $X_3, \dots$ of a group $A$, one defines inductively 
\begin{eqnarray*}
D_1 (X_1,X_2)& = & [X_1,X_2], \\
D_2 (X_1,X_2,X_3,X_4) & = & [[X_1,X_2],[X_3,X_4]],\\
\cdots & = & \cdots \\
D_k (X_1, \dots, X_{2^k}) & = & [D_{k-1} (X_1, \dots , X_{2^{k-1}}), D_{k-1}(X_{2^{k-1}+1}, \dots, X_{2^k})]. 
\end{eqnarray*}
Thus, in this situation, $J(A)=D_k(A,\dots,A)$. Furthermore, the arguments of Subsection~\ref{abelianization} are easily adapted to show that Lemma~\ref{joingroup} implies
\[
L_n[f] = \prod_{I_1\cup\dots\cup I_{2^k}=\set{n}} D_k\Bigl(\bigcap_{i\in I_1} K[f_i],\dots, \bigcap_{i\in I_{2^k}}K[f_i]\Bigr),
\]
for any $n\geq 0$ and any $n$-extension of groups $f$. Consequently, Theorem~\ref{Theorem-Hopf-BI} yields the following formulae:

\begin{theorem}
For any $n\geq 1$ and any $n$-presentation $f$ of a group $G$, an isomorphism
\begin{equation*}
H_{n+1} (G,\sol_k) \cong 
\frac{ D_k(f_{\varnothing}, \dots, f_{\varnothing}) \cap \bigcap_{i\in \set{n}}K[f_{i}]}{\prod_{I_1\cup\dots\cup I_{2^k}=\set{n}} D_k(\bigcap_{i\in I_1} K[f_i],\dots, \bigcap_{i\in I_{2^k}}K[f_i]) }
\end{equation*}
exists.\noproof 
\end{theorem}

\subsection{Some additional examples}
The reader will find it easy to further adapt the arguments used above to prove that Theorem~\ref{Theorem-Hopf-BI} yields similar formulae, for example in the variety of non-unital rings (where the role of commutator is played by the product of subrings) or in the variety of Lie algebras (the Lie bracket plays the role of commutator). The example of precrossed modules versus crossed modules, however, deserves special attention.

\subsection{Precrossed modules vs.\ crossed modules}
Recall that a \emph{precrossed} \emph{module} $(C,G,\del)$ is a group homomorphism $\del\colon C\to G$ equipped with a (left) group action of $G$ on $C$, such that
\[
\del(^gc) = g\del(c)g^{-1}
\]
for all $g\in G$ and $c\in C$. A morphism $f\colon {(C,G,\del)\to (D,H,\epsilon)}$ of precrossed modules is a pair of group homomorphisms $f_1\colon C\to D$ and $f_0\colon G\to H$ which preserve the action and are such that $\epsilon \comp f_1 =f_0\comp\del$. Let us write $\PXM$ for the category of precrossed modules and $\XM$ for the category of \emph{crossed} \emph{modules}, where this latter is the full subcategory of $\PXM$ whose objects $(C,G,\del)$ satisfy the \emph{Peiffer condition}
\[
^{\del(c)}c'=cc'c^{-1}
\]
for all $c, c'\in C$. We will sometimes abbreviate the notation $(C,G,\del)$ to $C$. 

It is well known that the category of precrossed modules is equivalent to a variety of $\Omega$-groups (see, e.g., \cite{LR} ~\cite{Loday} or ~\cite{Janelidze-Marki-Tholen}). Via this equivalence, $\XM$ correspond to a subvariety of $\PXM$. Hence, we can consider the Hopf formulae~\ref{Theorem-Hopf-BI} in the case where $\Ac=\PXM$ and $\Bc=\XM$.

Recall that a \emph{precrossed submodule} of a precrossed module $(C,G,\del)$ is a precrossed module  $(M,S,\mu)$ such that $M$ and $S$ are, respectively, subgroups of $C$ and $G$, and such that  the action of $S$ on $M$ is a restriction of the action of $G$ on $C$ and $\mu$ a restriction of $\del$ (in this case, we will write $\del$ instead of~$\mu$). $(M,S,\del)$ is a \emph{normal} precrossed submodule of $(C,G,\del)$ if, furthermore, $M$ and $S$ are, respectively, normal subgroups of $C$ and $G$, and, for all $c\in C$, $g\in G$, $m\in M$, $s\in S$, one has $^{g}m\in M$ and $^{s}cc^{-1}\in M$. This is exactly the case when $(M,S,\del)$ is the kernel of some morphism ${(C,G,\del)\to (D,H,\epsilon)}$. Note that the quotient ${(C,G,\del)\to(C,G,\del)/(K,S,\del)}$ is a degreewise quotient 
\[
(q_M,q_S)\colon {(C,G,\del)\to(C/M, G/S, \overline{\del})}.
\]
Here we have written $q_M$ and $q_S$ for the quotient homomorphisms ${C\to C/M}$ and ${G\to G/S}$, respectively. Note also that limits in $\PXM$ are degreewise limits in~$\Gp$.  

Let $(C,G,\del)$ be a precrossed module. The \emph{Peiffer commutator} of two precrossed submodules $(M,S,\del)$ and $(N,T,\del)$ of $(C,G,\del)$ is the normal subgroup of the group $M\cdot N$, generated by the \emph{Peiffer elements} $\langle m,n\rangle =mnm^{-1}(^{\del m}n)^{-1}$ and $\langle n,m\rangle =nmn^{-1}(^{\del n}m)^{-1}$, with $m\in M$ and $n\in N$. We will denote it by $\langle (M,S,\del),(N,T,\del)\rangle$ or simply by $\langle M,N\rangle$. Note that, for all $m\in M$ and $n\in N$, we have that 
\begin{eqnarray*}
\del \langle m,n\rangle &=& \del m \del n (\del m)^{-1} (\del {}^{\del m}n)^{-1}\\
&=& \del m \del n (\del m)^{-1} (\del m \del n (\del m)^{-1} )^{-1}\\
&=& 1.
\end{eqnarray*} 
Hence we may consider $\langle M,N\rangle$ as the precrossed submodule $(\langle M,N\rangle, 0,\del)$ of~$C$. 

The following simple identities, similar to~\ref{WH1} and~\ref{WH2}, are valid in any precrossed module $(C,G,\del)$, for any elements $a$, $a_1$, $a_2$, $b$, $b_1$, $b_2\in C$.
\begin{eqnarray}
\langle a,b_1 b_2\rangle & = & \langle a,b_1\rangle {}^{\del a}{b_1}\langle a,b_2\rangle ( ^{\del a}b_1)^{-1}  \label{PWH1}\\ 
\langle a_1a_2,b\rangle
& = & {a_1}\langle a_2,b\rangle a_1^{-1} \ \langle a_1, {}^{\del a_2}b \rangle \label{PWH2}
\end{eqnarray}
For a precrossed module $C$, let us denote by $\Delta_C=(1_{C},1_{C})\colon {C\to C\times C}$ the \emph{diagonal} of $C$: the smallest internal equivalence relation on $C$.

\begin{lemma}\label{joinmodule}
Suppose $(C,G,\del)$ is a precrossed module, $(M,S,\del)$ and $(N,T,\del)$ are normal precrossed submodules of $C$, and $(B,P,\del\times\del)$ is a precrossed submodule of $C\times C$. When $B$ is a either a normal precrossed submodule of $C\times C$ or $B=\Delta_K$, for some normal precrossed submodule $(K,R,\del)$ of $C$, then the following identity holds:
\[
\langle M\times_{\frac{M}{M\cap N}} M, B\rangle = \langle 0\times (M\cap N), B\rangle \cdot \langle\Delta_M, B\rangle.
\]
\end{lemma}
\begin{proof}
Note that in both cases only one inclusion is not entirely trivial. We prove the other.

Let us first assume that $B$ is a normal precrossed submodule of $C\times C$. Suppose that $m\in M$, $d\in M\cap N$ and $(b_1,b_2)\in B$, then by identity~\ref{PWH2}, 
\[
\langle (m,dm),(b_1,b_2) \rangle = \langle (1,d)(m,m) , (b_1,b_2) \rangle \in  \langle 0\times (M\cap N), B\rangle \cdot \langle\Delta_M, B\rangle.
\]
Similarly, by identity~\ref{PWH1}, $\langle (b_1,b_2), (m,dm)\rangle \in  \langle 0\times (M\cap N), B\rangle \cdot \langle\Delta_M, B\rangle$.

Let us now assume that $B=\Delta_K$, for some normal precrossed submodule $(K,R,\del)$ of $C$. Take $m\in M$, $d\in M\cap N$ and $k\in K$. It follows from \ref{PWH1} that 
\[
\langle (k,k),(m,dm)\rangle = \langle (k,k), (1,d)(m,m) \rangle\in  \langle 0\times (M\cap N), B\rangle \cdot \langle\Delta_M, B\rangle.
\]
Furthermore, 
\begin{eqnarray*}
\langle (m,dm),(k,k)\rangle&=& (\langle m,k\rangle ,\langle dm,k\rangle) \\
&=&  (\langle m,k\rangle , d\langle m,k\rangle d^{-1} \langle d, {}^{\del m}k\rangle) \\
&=& (\langle m,k\rangle,d\langle m,k\rangle d^{-1}) (1, \langle d, {}^{\del m}k\rangle).
\end{eqnarray*}
On the one hand,
\[
 (\langle m,k\rangle,d\langle m,k\rangle d^{-1}) = (1,d) \langle (m,m) , (k,k) \rangle (1,d)^{-1} \in \langle \Delta_M, B\rangle.
\]
On the other hand,  
\[
(1, \langle d, ^{\del m}k\rangle) = \langle (1,d),( ^{\del m}k, {}^{\del m}k) \rangle \in \langle 0\times (M\cap N), B\rangle.\qedhere
\]
\end{proof}

Using Lemma~\ref{joinmodule}, the reader will find it easy to adapt the arguments of Section~\ref{abelianization} to show that 
\[
L_n[f] = \Bigl( \prod_{I \subseteq \set{n}}\Bigl\langle \bigcap_{i\in I} K[f_i],   \bigcap_{i\notin I} K[f_i] \Bigr\rangle, 0 ,\del \Bigr),
\]
for any $n\geq 0$ and any $n$-extension of precrossed modules $f$. Consequently, Theorem~\ref{Theorem-Hopf-BI} becomes

\begin{theorem}\label{Thm-pxmod}
For any $n\geq 1$ and any $n$-presentation $f$ of a precrossed module $(C,G,\del)$, an isomorphism
\begin{equation*}
H_{n+1} ((C,G,\del),\xmod) \cong \frac{\langle f_{\varnothing}, f_{\varnothing}\rangle  \cap \bigcap_{i\in \set{n}}K[f_{i}]}{  \prod_{I \subseteq \set{n}}\Bigl\langle \bigcap_{i\in I} K[f_i],   \bigcap_{i\notin I} K[f_i] \Bigr\rangle  }
\end{equation*}
exists.\noproof 
\end{theorem}  

\pagebreak\printindex

%\bibliography{hopf}
%\bibliographystyle{amsplain}

%.bbl

\providecommand{\bysame}{\leavevmode\hbox to3em{\hrulefill}\thinspace}
\providecommand{\MR}{\relax\ifhmode\unskip\space\fi MR }
\providecommand{\MRhref}[2]{%
  \href{http://www.ams.org/mathscinet-getitem?mr=#1}{#2}
}
\providecommand{\href}[2]{#2}

\end{document}